\documentclass[twoside,11pt,abbrvbib,preprint]{article}

\usepackage{amsmath}
\usepackage{jmlr2e}

\makeatletter
\let\jmlr@origmaketitle\maketitle
\renewcommand{\maketitle}{\begingroup
  \long\def\@makefntext##1{\@setpar{\@@par\@tempdima \hsize
    \advance\@tempdima-15pt\parshape \@ne 15pt \@tempdima}\par
    \parindent 2em\noindent \hbox to \z@{\hss{\@thefnmark}\ \hfil}##1}%
  \jmlr@origmaketitle\endgroup}
\makeatother

\usepackage{bm}
\usepackage{xcolor}
\usepackage{algorithm}
\usepackage{algpseudocode}
\usepackage{tikz}
\usepackage{array}
\usepackage{enumitem}
\usepackage{outlines}

\newlist{plan}{enumerate}{3}
\setlist[plan,1]{label=\arabic*.}
\setlist[plan,2]{label=\arabic{plani}.\arabic*}
\setlist[plan,3]{label=\textbullet}

\newcolumntype{L}[1]{>{\raggedright\arraybackslash}p{#1}}

\hypersetup{colorlinks=true,linkcolor=blue,citecolor=blue,urlcolor=blue}
\usepackage[nosort]{cleveref}

\numberwithin{equation}{section}
\numberwithin{theorem}{section}

\usepackage{aliascnt}

\newaliascnt{lemma}{theorem}
\newtheorem{lemma}[lemma]{Lemma}
\aliascntresetthe{lemma}
\newaliascnt{corollary}{theorem}
\newtheorem{corollary}[corollary]{Corollary}
\aliascntresetthe{corollary}
\newaliascnt{remark}{theorem}
\newtheorem{remark}[remark]{Remark}
\aliascntresetthe{remark}
\usepackage{etoolbox}
\makeatletter
\AtBeginEnvironment{remark}{%
  \patchcmd{\@begintheorem}{\itshape}{\upshape}{}{}%
  \patchcmd{\NRorg@opargbegintheorem}{\itshape}{\upshape}{}{}}
\makeatother
\newaliascnt{proposition}{theorem}
\newtheorem{proposition}[proposition]{Proposition}
\aliascntresetthe{proposition}
\crefname{lemma}{lemma}{lemmas}
\Crefname{lemma}{Lemma}{Lemmas}
\crefname{corollary}{corollary}{corollaries}
\Crefname{corollary}{Corollary}{Corollaries}
\crefname{remark}{remark}{remarks}
\Crefname{remark}{Remark}{Remarks}
\crefname{proposition}{proposition}{propositions}
\Crefname{proposition}{Proposition}{Propositions}
\crefname{equation}{eq.}{eqns.}
\Crefname{equation}{Eq.}{Eqns.}
\newtheorem{assumption}{Assumption}
\newtheorem{fact}{Fact}
\crefname{fact}{fact}{facts}
\Crefname{fact}{Fact}{Facts}

\newcommand{\R}{\mathbb{R}}

\newcommand{\E}{\mathbb{E}}

\newcommand{\prox}{\operatorname{prox}}

\newcommand{\vm}{\mathbf{m}}
\newcommand{\vs}{\mathbf{s}}
\newcommand{\vx}{\mathbf{x}}
\newcommand{\vy}{\mathbf{y}}
\newcommand{\vz}{\mathbf{z}}
\newcommand{\vu}{\mathbf{u}}

\newcommand{\vq}{\mathbf{q}}
\newcommand{\vw}{\mathbf{w}}
\newcommand{\va}{\mathbf{a}}
\newcommand{\vb}{\mathbf{b}}

\newcommand{\grad}{\mathbf{\nabla}}
\newcommand{\sgrad}{\widetilde{\mathbf{\nabla}}}

\crefname{assumption}{assumption}{assumptions}
\Crefname{assumption}{Assumption}{Assumptions}

\usepackage{lastpage}

\ShortHeadings{Last-Iterate Convergence of SGD for Unbounded Problems}{Vaidyan, Alacaoglu, Friedlander}
\firstpageno{1}

\begin{document}

\title{Last-Iterate Convergence of Stochastic Gradient Methods for Unbounded Problems}
  
\author{\name Kevin Kurian Thomas Vaidyan \email kevinktv@cs.ubc.ca \\
       \addr Department of Computer Science\\
       University of British Columbia
       \AND
       \name Ahmet Alacaoglu\thanks{Joint last authors.} \email ahmet.alacaoglu@ubc.ca \\
       \addr Department of Mathematics\\
       University of British Columbia
              \AND
       \name Michael P. Friedlander\footnotemark[1] \email michael.friedlander@ubc.ca \\
       \addr Department of Computer Science and Mathematics\\
       University of British Columbia}

\maketitle

\begin{abstract}%
Last-iterate convergence guarantees for stochastic gradient methods have attracted significant attention as practical implementations typically return the last iterate instead of an average of the iterates. In this paper, we study two approaches for analyzing the last-iterate convergence of proximal and projected stochastic gradient methods without bounded domain or bounded variance assumptions, where black-box variational inference serves as a running example. Under two related assumptions relaxing the bounded variance, we show near-optimal anytime rate $\widetilde{\mathcal{O}}(1/\sqrt{T})$ for the last iterate. Moreover, when the horizon $T$ is known, we remove the logarithmic factor and obtain the optimal order $\mathcal{O}(1/\sqrt{T})$, matching the order known under bounded variance or bounded domain assumptions.
\end{abstract}


\section{Introduction}\label{sec:Intro}
Stochastic gradient descent (SGD) and its variants are widely used in training large-scale machine learning models. In this paper, we study two standard stochastic first-order methods: proximal stochastic gradient descent (proximal SGD)
for regularized problems of the form
    \begin{equation}\label{eq: regularized_problem}
    \min_{\vx\in\R^d} \big\{h(\vx) := f(\vx)+g(\vx)\big\},
    \end{equation}
    and projected stochastic gradient descent (projected SGD) 
    for convex constrained optimization problems of the form
    \begin{equation}\label{eq: constr_problem}
    \min_{\vx\in \mathcal{C}} f(\vx).
    \end{equation}
Most classical analyses of SGD provide guarantees on the suboptimality of averaged iterates or their weighted combinations  \citep{nemirovski2009robust,lan2020first}. That is, it is common to obtain rates for $f(\bar \vx_T)-f(\vx^*)$ where $\bar \vx_T = \frac{1}{T} \sum_{t=1}^T \vx_t$.
However, practitioners typically use the last iterate as the output, rather than the averaged iterate \citep{shamir2013stochastic,jain2019making}. For example, in the $\ell_1$-regularized least-squares problem
\begin{equation}\label{eq: slo4}
    \min_{\vx \in \R^d} \frac{1}{n} \sum_{i=1}^n (\va_i^\top \vx - b_i)^2 + \lambda \|\vx\|_1,
\end{equation}
the $\ell_1$ regularizer promotes sparsity in the solution. However, averaging can destroy this sparsity; hence, guarantees for the last iterate $\vx_T$ are of interest \citep[Sec. 1]{chen2012optimal}.

A second challenge concerns the assumptions made on the stochastic oracle. Many analyses require strong assumptions such as bounded second moments or bounded variance of the stochastic gradients $\sgrad f(\vx)$ \citep{robbins1951stochastic,nemirovski2009robust,lan2020first}. That is, these works assume
    \begin{equation}\label{eq: bdd_grad_var}
    \mathbb{E} \| \sgrad f(\vx)\|^2\leq \sigma^2 \text{~~~or~~~} \mathbb{E} \| \sgrad f(\vx)-\nabla f(\vx)\|^2\leq \sigma^2.
    \end{equation}
Both of these assumptions are widely accepted to be rather restrictive and may not hold over unbounded domains even for simple problems such as \eqref{eq: slo4}. 
In the context of \eqref{eq: slo4}, the stochastic gradient is given by $\sgrad f(\vx) = 2\va_i(\langle \va_i, \vx \rangle - b_i)$, which grows with $\|\vx\|$ and is unbounded.

These observations motivate the central question of this paper:
\emph{Can SGD (and its proximal and projected variants) attain
last-iterate convergence guarantees without assuming uniformly bounded
variance of the stochastic gradient oracle?}

We answer this question under two complementary conditions: component-wise convexity and smoothness \eqref{eq: cw_cond} (\Cref{sec: componentwise_analysis}) and the Blum-Gladyshev \eqref{eq: BG} growth condition (\Cref{sec: bg_analysis}). 
Under both conditions, we establish an anytime last-iterate convergence rate of $\mathcal{O}(\ln T/\sqrt{T})$ and a log-free rate of $\mathcal{O}(1/\sqrt{T})$ when the horizon is known. Our analysis under \eqref{eq: BG} accommodates convex objectives that need not be smooth or Lipschitz continuous and yields applications to black-box variational inference. Under \eqref{eq: cw_cond}, we extend the existing proximal SGD guarantees of \citet{vaidyan2026convergencerate} to variable step sizes and improve the convergence rate for horizon-dependent step sizes.

\paragraph{A Brief Discussion on Prior Work.}
Classical analyses established convergence rates for averaged iterates under bounded domain and bounded second-moment assumptions \citep{nemirovski2009robust}. Subsequent works relaxed these assumptions, allowing unbounded domains and growing stochastic gradients \citep{grimmer2019convergence,khaled2023unified}. 
A parallel line of research studied last-iterate convergence, under bounded second-moment or bounded variance assumptions \citep{shamir2013stochastic,harvey2019tight,jain2019making,orabona2020last,liu2024revisitinglastiterateconvergencestochastic}. 
More recently, several works have obtained last-iterate convergence guarantees under weaker assumptions on the variance \citep{garrigos2025last,kim2025nearly,attia2025fast,alacaoglu2025towards,vaidyan2026convergencerate}. 
More detailed comparisons appear in the respective sections in the sequel.

\paragraph{Component-Wise Convexity and Smoothness.}    The first approach to weaken \eqref{eq: bdd_grad_var} is to assume component-wise convexity and smoothness of the objective \citep{garrigos2025last,attia2025fast}. In this setting, the objective decomposes as
    \begin{equation}\label{eq: stoc_func}
    f(\vx) = \E_{i\sim P}[f_i(\vx)],
    \end{equation}
for some unknown distribution $P$, from which we get i.i.d. samples. In particular, we assume that
\begin{equation}\tag{CW}\label{eq: cw_cond}
\text{each component function~} f_i  \text{~is convex and $L$-smooth.}
\end{equation} 
    The most standard setting where this assumption holds is the least squares problem with a regularizer; see \eqref{eq: slo4}, where each component $f_i(\vx) = (\va_i^\top\vx - b_i)^2$ is convex and $(\max_i 2\|\va_i\|^2)$-smooth.
    \citet{garrigos2025last} and \citet{attia2025fast} used this assumption along with the powerful technique introduced by \citet{zamani2025exact} to obtain last-iterate guarantees for \eqref{eq: regularized_problem} with $g \equiv 0$. A classical consequence of the component-wise convexity and smoothness assumption is the bound
    \begin{equation}\label{eq: ab_type}
    \mathbb{E}\|\grad f_i(\vx)\|^2 \leq 2\mathbb{E}\|\grad f_i(\vx^*)\|^2 + 4L(f(\vx)-f(\vx^*)),
    \end{equation}
    where $\vx^*$ is a minimizer of $f$. This bound is a special case of the so-called ABC assumption \citep{polyak1973pseudogradient,khaled2022better,khaled2023unified}. 
    
    Using this property, \citet{garrigos2025last} and \citet{attia2025fast} show that the last iterate of SGD converges at a rate of $\mathcal{O}(\ln T/\sqrt{T})$. The analysis developed in these works
    is restricted to horizon-dependent step sizes and does not handle regularized problems such as \eqref{eq: regularized_problem} where the proximal method is needed. The work by \citet{vaidyan2026convergencerate} handles the proximal setting but is still restricted to horizon-dependent step sizes and suffers from a logarithmic factor in the rate, as do earlier works in this line.

\paragraph{Blum-Gladyshev Variance.}    The second approach to weaken \eqref{eq: bdd_grad_var} is the Blum-Gladyshev (BG) condition \citep{blum1954approximation, gladyshev1965stochastic} given as
    \begin{equation}\tag{B-G}\label{eq: BG}
    \mathbb{E}\|\sgrad f(\vx)\|^2 \leq B^2 \|\vx-\vx_0\|^2 + G^2,
    \end{equation}
    where $\vx_0$ is the initial iterate of the algorithm. Note that this assumption is also sometimes presented with a solution $\vx^*$ instead of $\vx_0$ on the right-hand side; the two forms are interchangeable (see \Cref{fact: bg centering}). 
    
    Unlike \eqref{eq: cw_cond}, the Blum-Gladyshev condition \eqref{eq: BG} does not require the individual stochastic functions $f_i$ to be smooth, Lipschitz continuous, or convex. We provide more concrete examples to illustrate this point in \Cref{sec: examples}. 
    
    Moreover, the BG condition is strictly weaker than the ABC assumption.
    It was shown, for example, by \citet{alacaoglu2025towards}, that \eqref{eq: ab_type} (which itself is implied by \ref{eq: cw_cond}) implies \eqref{eq: BG}, when $f$ is $L$-smooth, by straightforward estimations. 

    However, the converse need not be true. \citet[Example 1.1]{alacaoglu2025towards} show that, for a simple quadratic example over an unbounded domain with a standard stochastic oracle, \eqref{eq: BG} is satisfied, but \eqref{eq: ab_type} (and, consequently, \ref{eq: cw_cond}) fails to hold.
This establishes that \eqref{eq: BG} is a weaker condition than \eqref{eq: cw_cond} even for smooth problems.

Additionally, without the expectation, the \eqref{eq: BG} assumption also describes a class of functions $f$ that are \emph{neither smooth nor Lipschitz}.
For this class, we obtain last-iterate rates for the classical subgradient method (see \Cref{rmk: deterministic subgradient descent}), which are new to the best of our knowledge.
Existing results for this class either give rates on the ergodic average \citep{lu2019relative,grimmer2019convergence,li2025revisiting} or rely on transformations of the original problem to design specialized algorithms that depart from the simple subgradient method \citep{renegar2016efficient,grimmer2018radial}. 

\paragraph{Application: Black-Box Variational Inference.} 
Interestingly, the \eqref{eq: BG} assumption has recently been shown to be relevant for black-box variational inference (BBVI) \citep{kim2025nearly,labarriere2026sgd,domke2023provable}. 
Variational inference is a classical technique in statistics and machine learning, widely used in applications including healthcare \citep{ranganath2016deep, lopez2018deep}. Its goal is to approximate a potentially intractable posterior distribution by a distribution from a simpler family, such as Gaussians, by minimizing the Kullback–Leibler divergence between the two.

When the target distribution is complex, closed-form updates for variational inference are often unavailable. However, accessing stochastic gradients is still feasible, making stochastic optimization-based approaches, such as BBVI \citep{ranganath2014black}, popular for this task. 
In BBVI, the second moment of the standard gradient estimator is not uniformly bounded, so \eqref{eq: bdd_grad_var} does not hold. However, under smoothness assumptions on the log-posterior, it satisfies the \eqref{eq: BG} condition \citep{domke2023provable}. See \Cref{subsec: bbvi} for technical details relevant to our exposition. 

\paragraph{Comparison of the Two Approaches.}
The assumptions \eqref{eq: BG} and \eqref{eq: cw_cond} are complementary. As discussed earlier, the \eqref{eq: BG} assumption is weaker and covers a broader class of problems not handled by \eqref{eq: cw_cond}.  However, under \eqref{eq: cw_cond}, the standard step size $\tau_t=\tau/\sqrt{t+1}$ attains the optimal rate up to a logarithmic factor, which our analysis under \eqref{eq: BG} does not support; see \Cref{remark: oracle_comparison}. Hence, there is a trade-off between the two assumptions: \eqref{eq: BG} is weaker but yields suboptimal order in the rate with the standard step size of $\tau_t=\tau/\sqrt{t+1}$, while \eqref{eq: cw_cond} is stronger but yields the optimal order in a wider variety of settings. See also Table \ref{tab:assumption_comparison}.

\paragraph{Summary of Contributions.}
\begin{itemize}
    \item \textbf{Last-iterate rates under \eqref{eq: BG}.} We establish last-iterate convergence rates for proximal SGD when $f$ is convex but not necessarily smooth or Lipschitz. Specifically, we obtain an anytime rate of $\mathcal{O}(\ln T/\sqrt{T})$ and, when the horizon is known, a log-free rate of $\mathcal{O}(1/\sqrt{T})$ using the step-size scheme of \citet{liu2024revisitinglastiterateconvergencestochastic} (see \Cref{theorem: proxsgd_bg_nonsmooth} for more details). Our results also apply to projected SGD. To our knowledge, these are the first anytime and the first log-free last-iterate rates for SGD without assuming bounded variance. As a special case, we obtain last-iterate rates for the classical subgradient method on convex functions that are not Lipschitz.
    \item \textbf{Black-box variational inference.} Using the gradient bounds of \citet{domke2023provable}, our results give anytime last-iterate guarantees for proximal and projected SGD for BBVI (\Cref{subsec: bbvi}). Existing guarantees hold for weighted averages of the iterates with horizon-dependent step sizes. Our results answer an open question by \citet[Section~2.4]{kim2025nearly} regarding anytime guarantees of SGD for BBVI without strong convexity.
    \item \textbf{Variable step sizes under \eqref{eq: cw_cond}.} We extend the proximal SGD analysis of \citet{vaidyan2026convergencerate}, which is restricted to constant horizon-dependent step sizes, to handle non-increasing variable step sizes. We establish an anytime convergence rate of $\mathcal{O}(\ln T/\sqrt{T})$ and a log-free rate of $\mathcal{O}(1/\sqrt{T})$ when the horizon is known (see \Cref{theorem: convergenceratesforproximalsgd_componentwise} for more details). Projected SGD follows as a special case. These results answer an open question posed by \cite[Remark 3.6, Section 5]{garrigos2025last}, in the unconstrained case.
\end{itemize}

\begin{table}[!htb]
\centering
\small
{\begin{tabular}{L{0.24\textwidth} | L{0.30\textwidth} L{0.37\textwidth}}
\hline
\rule{0pt}{5mm}\textbf{Property}
& \textbf{\eqref{eq: cw_cond} Assumption}
& \textbf{\eqref{eq: BG} Assumption} \\[2mm]
\hline
Allows nonconvex $f_i$ 
& $\times$
& $\checkmark$ \\[3mm]
Allows nonsmooth $f_i$
& $\times$
& $\checkmark$ \\[3mm]
Control mechanism
& ABC-type inequality \eqref{eq: ab_type}
& Iterates bounded in expectation \\[3mm]
Step size regimes
& $\tau_t \le \frac{1}{4L\sqrt{t+1}}$
& $\sum_{t=1}^T \tau_t^2=\mathcal{O}(1)$
\end{tabular}}
\caption{Comparison of the assumptions considered in this paper.\label{tab:assumption_comparison}}
\end{table}

\paragraph{Analysis Roadmap.} We analyze SGD-based methods under both \eqref{eq: BG} and \eqref{eq: cw_cond}. In both cases, the analysis unfolds in two steps:
\begin{enumerate}
  \item We first establish a control on the second moment of the stochastic gradients along the trajectory of the algorithm, that is, on the quantity $\mathbb{E}\|\widetilde{\nabla} f(\vx_t)\|^2$. Under \eqref{eq: cw_cond}, this control is obtained through the ABC condition in \eqref{eq: ab_type}. Under \eqref{eq: BG}, we show that as long as the sum of the squared step sizes does not grow too quickly, the iterates of SGD remain bounded in expectation, which implies that the second moment of the stochastic gradients is bounded along the trajectory of the algorithm.
  \item With this control in place, we adapt the reference point argument of \citet{zamani2025exact} and its generalization in \citet{liu2024revisitinglastiterateconvergencestochastic}: we prove a one-step progress bound against an arbitrary reference point, then choose these points as convex combinations of the minimizer and past iterates so that summing the bounds gets us the last-iterate guarantee.
\end{enumerate}

\paragraph{Relation to \citet{vaidyan2026convergencerate}.}
This paper extends our conference paper \citep{vaidyan2026convergencerate}, which analyzed proximal SGD under \eqref{eq: cw_cond} with a constant, horizon-dependent step size and obtained an $\mathcal{O}(\ln T/\sqrt{T})$ last-iterate rate. The present paper differs from it in three ways (see also Table \ref{tab:related_work_comparison}). 
\begin{itemize}
    \item The analysis under \eqref{eq: BG} in \Cref{sec: bg_analysis} is new. It controls the second moment of the stochastic gradients through boundedness of the iterates in expectation, rather than through the ABC-type inequality \eqref{eq: ab_type} used in \citet{vaidyan2026convergencerate}. As a result, it applies to proximal SGD on convex objectives that need not be smooth or Lipschitz.
    \item The applications in \Cref{sec: examples} are new. We show that \eqref{eq: BG} covers problems outside the scope of \eqref{eq: cw_cond}, such as minimizing the maximum of convex quadratics (\Cref{subsec: nonsmooth_example}), and we use it to obtain last-iterate guarantees for SGD in black-box variational inference (\Cref{subsec: bbvi}). This application was not covered in \citet{vaidyan2026convergencerate}, because the \eqref{eq: cw_cond} assumption does not hold for BBVI.
    \item In \Cref{sec: componentwise_analysis} we return to the \eqref{eq: cw_cond} setting of \citet{vaidyan2026convergencerate} and extend its analysis from constant to non-increasing variable step sizes. The key technical step, \Cref{theorem: proxsgd_decreasingstepsize}, generalizes the reference-point construction of \citet[Lemma A.6]{vaidyan2026convergencerate}; this gives an anytime $\mathcal{O}(\ln T/\sqrt{T})$ rate with $\tau_t = \tau/\sqrt{t+1}$, and a log-free $\mathcal{O}(1/\sqrt{T})$ rate when the horizon is known (\Cref{theorem: convergenceratesforproximalsgd_componentwise}). This improves on the guarantee in \citet[Theorem 3.1]{vaidyan2026convergencerate}, which required constant step sizes and included a log term.
\end{itemize}

\begin{table}[!tb]
\centering
\small
{\begin{tabular}{L{0.2\textwidth} L{0.17\textwidth} L{0.23\textwidth} L{0.13\textwidth} L{0.13\textwidth}}
\hline
\textbf{Reference}
& \textbf{Constrained / Regularized}
& \textbf{Assumptions}
& \textbf{Anytime rate}
& \textbf{Horizon-dependent rate} \\
\hline
\citet{shamir2013stochastic}
& \checkmark
& Bdd gradients/domain
& $\mathcal{O}\!\left(\frac{\ln T}{\sqrt{T}}\right)$
& -- \\[2mm]
\citet{jain2019making}
& \checkmark
& Bdd gradients/domain
& --
& $\mathcal{O}\!\left(\frac{1}{\sqrt{T}}\right)$ \\[4mm]
\citet{liu2024revisitinglastiterateconvergencestochastic}
& \checkmark
& Bdd variance
& $\mathcal{O}\!\left(\frac{\ln T}{\sqrt{T}}\right)$
& $\mathcal{O}\!\left(\frac{1}{\sqrt{T}}\right)$ \\[2mm]
\citet{garrigos2025last,attia2025fast}
& $\times$
& \eqref{eq: cw_cond}
& --
& $\mathcal{O}\!\left(\frac{\ln T}{\sqrt{T}}\right)$ \\[2mm]
\citet{vaidyan2026convergencerate}
& \checkmark
& \eqref{eq: cw_cond}
& --
& $\mathcal{O}\!\left(\frac{\ln T}{\sqrt{T}}\right)$ \\[2mm]
\hline
This work
& \checkmark
& \eqref{eq: BG}
& $\mathcal{O}\!\left(\frac{\ln T}{\sqrt{T}}\right)$
& $\mathcal{O}\!\left(\frac{1}{\sqrt{T}}\right)$ \\[2mm]
This work
& \checkmark
& \eqref{eq: cw_cond}
& $\mathcal{O}\!\left(\frac{\ln T}{\sqrt{T}}\right)$
& $\mathcal{O}\!\left(\frac{1}{\sqrt{T}}\right)$ \\
\hline
\end{tabular}}
\smallskip
\caption{Comparison with existing last-iterate convergence results for SGD variants.
\label{tab:related_work_comparison}}
\end{table}

\section{Problem Setup and Assumptions}\label{sec: setup}
In this section, we introduce the two stochastic first-order methods studied in the paper and fix the notation used throughout.
Given iterates $(\vx_t)$, let $\mathcal{F}_t := \sigma(\vx_0, \dots, \vx_t)$ denote the $\sigma$-algebra generated by $\vx_0, \dots, \vx_t$, and write $\E_t[\cdot] := \E[\cdot \mid \mathcal{F}_t]$ for the conditional expectation with respect to $\mathcal{F}_t$. We use $\langle \cdot, \cdot \rangle$ for the Euclidean inner product and $\|\cdot\|$ for the Euclidean norm. We use $\partial f(\vx)$ to denote the subdifferential of a convex function $f$ at $\vx$.

\paragraph{Proximal Stochastic Gradient Descent.}
For the additively composite problem given in \eqref{eq: regularized_problem},
we assume that $f\colon\R^d \to \R$ is convex, but not necessarily smooth, and $g\colon\R^d \to \R \cup \{+\infty\}$ is proper, convex and closed. We denote $\vx^* \in \arg \min_{\vx \in \R^d} h(\vx)$ as an optimal solution and $h^* = h(\vx^*)$. The proximal operator of $g$ is defined by
    \begin{equation*}
    \prox_{g}(\vy)
    =
    \arg\min_{\vx \in \R^d}
    \left\{
    g(\vx)
    +
    \frac{1}{2}\|\vx-\vy\|^2
    \right\}.
    \end{equation*}
    Given $\vx_0\in \operatorname{dom}g$, the proximal SGD update is given for $t\geq 0$ as
    \begin{equation}\tag{proxSGD}\label{eq: prox_sgd_update}
    \vx_{t+1}
    =
    \prox_{\tau_t g}\left(\vx_t - \tau_t \sgrad f(\vx_t)\right),
    \end{equation}
    where $\sgrad f(\vx_t)$ denotes an unbiased stochastic subgradient of $f$ satisfying
    \begin{equation*}
    \E_t[\sgrad f(\vx_t)] \in \partial f(\vx_t).
    \end{equation*}

\paragraph{Projected Stochastic Gradient Descent.}\label{sec: Proj SGD summary}
The constrained problem \eqref{eq: constr_problem} is the special case of \eqref{eq: regularized_problem} with $g = \delta_{\mathcal{C}}$, where $\delta_{\mathcal{C}}$ is the indicator function of a nonempty, closed and convex set $\mathcal{C}$, defined as $\delta_{\mathcal{C}}(\vx) = 0$ if $\vx \in \mathcal{C}$ and $\delta_{\mathcal{C}}(\vx) = +\infty$ otherwise. 
In this case, $\prox_{\tau_t g}$ is the Euclidean projection onto $\mathcal{C}$, so \eqref{eq: prox_sgd_update} reduces to projected SGD, and all of our results for proximal SGD apply to it.
    \paragraph{Assumptions.} We now state the standard assumptions under which we analyze proximal SGD.
    \begin{assumption}\label{assumption: prox}
        The problem \eqref{eq: regularized_problem} satisfies the following:
        \begin{enumerate}
            \item The solution set of \eqref{eq: regularized_problem} is nonempty.
            \item The function $g\colon\R^d\to\R\cup\{+\infty\}$ is proper, closed and convex.
            \item The function $f\colon\R^d\to\R$ is convex.
            \item The oracle is unbiased: $\E_t[\sgrad f(\vx_t)] \in \partial f(\vx_t)$ for all $t\geq 0$.
        \end{enumerate}
    \end{assumption}
    Note that \Cref{assumption: prox} does not require $f$ to be differentiable. For projected SGD, that is, with $g = \delta_{\mathcal{C}}$, the second item amounts to $\mathcal{C}$ being nonempty, closed and convex.
    \Cref{assumption: prox} is standard in the literature on projected SGD \citep{nemirovski2009robust,shamir2013stochastic,grimmer2019convergence} and proximal SGD \citep{khaled2023unified,rosasco2020convergence,vaidyan2026convergencerate}.
    For the stochastic oracle, we consider the two models discussed in \Cref{sec:Intro}.
    \begin{assumption}\label{assumption: oracle}
        The stochastic oracle satisfies one of the following:
        \begin{enumerate}[label=(\alph*), ref=\theassumption(\alph*)]
            \item\label[assumption]{assumption: oracle BG} There exist $B, G \geq 0$ such that \eqref{eq: BG} holds for all $\vx\in\operatorname{dom}g$. That is, along the iterates, it holds that
            \[
                \E_t\|\sgrad f(\vx_t)\|^2 \leq B^2\|\vx_t-\vx_0\|^2 + G^2 .
            \]
            \item\label[assumption]{assumption: oracle CW} The following hold:
            \begin{enumerate}[label=(\roman*), ref=\theassumption(\alph{enumi})(\roman*)]
                \item The function $f$ decomposes as \eqref{eq: stoc_func}, where every $f_i$ satisfies \eqref{eq: cw_cond}.
                \item The oracle is $\sgrad f(\vx_t) = \grad f_{i_t}(\vx_t)$, where $i_t\sim P$ is sampled independently of $\mathcal{F}_t$.
                \item $\E_{i\sim P}\|\grad f_i(\vx^*)\|^2 \leq \sigma_*^2  < \infty$ for some solution $\vx^*$.
            \end{enumerate}
        \end{enumerate}
    \end{assumption}
    In the subsequent sections, we will be combining \Cref{assumption: prox} with an assumption on the oracle (\Cref{assumption: oracle BG} or \Cref{assumption: oracle CW}) to analyze the convergence of proximal SGD.
    Under \Cref{assumption: oracle CW}, $f$ is $L$-smooth since every $f_i$ is $L$-smooth, and $\E_t[\sgrad f(\vx_t)] = \grad f(\vx_t)$.

\section{Convergence Analysis under B-G}\label{sec: bg_analysis}
We now proceed to analyze the last-iterate convergence of proximal SGD under the BG condition \eqref{eq: BG}. Because of its applicability to BBVI, SGD under this assumption has been studied in several recent works \citep{domke2023provable,labarriere2026sgd,kim2025nearly}. The analyses in these works set up the recursion
\begin{equation}\label{eq: non_telescop}
2\tau_t \E\left[ h(\vx_{t+1}) - h^* \right] \leq (1+2\tau_t^2B^2)\E\|\vx_{t}-\vx^*\|^2 -\E\|\vx_{t+1}-\vx^*\|^2+2\tau_t^2G^2,
\end{equation}
where $\tau_t = \tau$ is fixed. In this recursion, the terms involving $\E\|\vx_t-\vx^*\|^2$ do not telescope due to the additional $1+ 2\tau^2B^2$ factor. These works then handle this issue by weighting the inequalities with $\alpha_{t+1} = \alpha_t/(1 + 2\tau^2B^2)$ which then leads to a telescoping sum on the right-hand side and yields a guarantee for the weighted average of the iterates. However, this approach does not accommodate variable step sizes, and it gives neither anytime nor last-iterate guarantees, both of which our analysis provides. 

\citet{li2025revisiting} studied the subgradient method without assuming Lipschitz continuity, using a strategy similar to our analysis roadmap. This work used normalized step sizes $\tau_t = \beta_t/\|\vs_t\|$ with $\vs_t\in\partial f(\vx_t)$ and $\sum_{k=0}^\infty \beta_k^2 < \infty$, which allows bounding the iterates and, consequently, the subgradient norms along the trajectory. With this bound in place, the standard analysis for bounded subgradients applies. Unlike our analysis, this approach requires the step size to depend on the norm of the sampled subgradient, which does not extend directly to stochastic oracles. As a result, this work analyzed the stochastic subgradient method only under the interpolation condition. In addition, this work derived rates for the average of the iterates, whereas our focus is on last-iterate rates.

Throughout this section, we require \Cref{assumption: prox,assumption: oracle BG}.
We follow the two-step strategy outlined in the analysis roadmap in \Cref{sec:Intro}. 
For the analysis, it is convenient to recenter the BG condition \eqref{eq: BG} at a solution $\vx^*$. We apply \Cref{fact: bg centering} with $\va = \vx_0$ and $\vb = \vx^*$ to obtain
\begin{equation}\label{eq: BG recentered}
\E_t\|\sgrad f(\vx_t)\|^2 \le 2B^2\|\vx_t-\vx^*\|^2 + \widehat G^2,
\qquad \widehat G^2 := 2B^2\|\vx_0-\vx^*\|^2 + G^2 .
\end{equation}
In order to establish boundedness of the iterates and the second moment of the stochastic subgradients, we first need the following lemma which analyzes one step of the algorithm.
\begin{lemma}\label{lemma: onestepid_proximal_nonsmooth}
    Let \Cref{assumption: prox} hold. For all $t\geq 0$, any $\mathcal{F}_t$-measurable $\vz_t$, and any $\vq_{t+1} \in \partial f(\vx_{t+1})$, we have
    \[\begin{aligned}
    2 \tau_t \E [h(\vx_{t+1}) - h(\vz_t)] \leq \E\|\vx_t - \vz_t\|^2 - \E \|\vx_{t+1} - \vz_t\|^2 + \tau_t^2 \E \|\vq_{t+1} - \sgrad f(\vx_t)\|^2.
    \end{aligned}\]
\end{lemma}
The proof of \Cref{lemma: onestepid_proximal_nonsmooth} is given in \Cref{app: onestepid_proximal_nonsmooth_proof}.
We now use the above lemma to establish the expected boundedness of the iterates and the second moment of the stochastic gradients along the trajectory of the algorithm.
\begin{proposition}[Expected boundedness]\label{prop: BG boundedness}
    Let \Cref{assumption: prox,assumption: oracle BG} hold, and set $\tau_t \leq \frac{1}{2\sqrt{2}B}$.
     If for all $t\geq 0$, the iterates satisfy
    \begin{equation}\label{eq: distance recursion}
        \E\|\vx_{t+1}-\vx^*\|^2 \leq \E\|\vx_t-\vx^*\|^2 + \tau_t^2\,\E\|\vq_{t+1} - \sgrad f(\vx_t)\|^2,
    \end{equation}
    where $\vq_{t+1} := \E_{t+1}[\sgrad f(\vx_{t+1})] \in \partial f(\vx_{t+1})$, then, we have
    \begin{equation*}
        \E \|\vx_{t+1} - \vx^*\|^2 \leq \exp\left( 16B^2 P_t \right) \left( \|\vx_0 - \vx^*\|^2 + 8 P_t \widehat{G}^2 \right) \text{~~and~~} \E\|\sgrad f(\vx_t)\|^2 \leq K_{t-1},
    \end{equation*}
    where $P_t := \sum_{i=0}^t \tau_i^2$ with $P_{-1} := 0$
    and $K_t := 2B^2\exp(16B^2P_t)\left(\|\vx_0-\vx^*\|^2 + 8P_t\widehat G^2\right) + \widehat G^2$.
\end{proposition}
See \Cref{app: BG_prop1} for the proof of \Cref{prop: BG boundedness}.
We emphasize that, with an appropriate choice of step size, the factor $\exp(16B^2P_t)$ is bounded by an absolute constant that grows with neither $T$ nor $B$. To see this, consider the step sizes in \Cref{theorem: proxsgd_bg_nonsmooth}(i)-(ii), where the quantity $P_t$ stays bounded by a constant multiple of $\tau^2$ (\Cref{proxboundsforstep2,proxboundsforstep4}(i)). Since $\tau$ is chosen proportional to $1/B$, the $B^2$ in the exponent cancels, so $B^2P_t$ is bounded by an absolute constant. 

For the third step size choice, that is, $\tau_t = \tau/\sqrt{t+1}$, the quantity $P_t$ grows logarithmically (\Cref{proxboundsforstep1}(i)), so this exponential factor grows polynomially in $T$, which leads to the suboptimal order in \Cref{theorem: proxsgd_bg_nonsmooth}(iii).

Applying \Cref{lemma: onestepid_proximal_nonsmooth} with $\vz_t = \vx^*$ and $\vq_{t+1} = \E_{t+1}[\sgrad f(\vx_{t+1})]$, which lies in $\partial f(\vx_{t+1})$ by \Cref{assumption: prox}, and dropping the nonnegative left-hand side shows that \eqref{eq: prox_sgd_update} satisfies \eqref{eq: distance recursion}. As a result, if $\tau_t \leq \frac{1}{2\sqrt{2}B}$, \Cref{prop: BG boundedness} applies, giving us the estimate $\E\|\sgrad f(\vx_t)\|^2 \leq K_{t-1}$ for all $0 \leq t$.

Thus, we have established control over the second moment, which is the first ingredient in the analysis roadmap. We now proceed to the second step in our roadmap. We obtain a one-step progress bound for proximal SGD under the BG noise model.
\begin{lemma}[One-step of proximal SGD under \eqref{eq: BG}]\label{proxobjectivecombBG nonsmooth}
Under \Cref{assumption: prox,assumption: oracle BG}, set $\tau_t \leq \frac{1}{2\sqrt{2}B}$.
    Then for all $t\geq 0$ and every $\mathcal{F}_t$-measurable $\vz_t$, \eqref{eq: prox_sgd_update} satisfies
    \begin{equation*}
        \E [h(\vx_{t+1}) - h(\vz_t)] \leq \frac{1}{2\tau_t} \E[\|\vx_t - \vz_t\|^2 - \|\vx_{t+1} - \vz_t\|^2] + 2\tau_t K_t,
    \end{equation*}
    where $K_t$ is defined in \Cref{prop: BG boundedness}.
\end{lemma}
The proof of \Cref{proxobjectivecombBG nonsmooth} is given in \Cref{app: BG_lem_nonsmooth}.

To convert this one-step bound into a last-iterate convergence rate, we use the reference-point technique of \citet{zamani2025exact}. Specifically, we take each $\vz_t$ to be a convex combination of $\vx^*$ and the past iterates, with weights chosen so that, after summing the one-step bounds, every objective gap except the last one cancels. We can apply this due to our two-step analysis: in the first step (see \Cref{prop: BG boundedness}), we established the bound $K_t$ on the second moment of the stochastic subgradients observed by the algorithm. Hence, we have reduced our setting to one similar to the \emph{bounded subgradient} settings of earlier work, by restricting to the trajectory of the algorithm. \Cref{theorem: proxsgd_decreasingstepsize} extends this technique to the proximal setting, building on the analysis of \citet{liu2024revisitinglastiterateconvergencestochastic}, which assumes bounded variance.

We state \Cref{theorem: proxsgd_decreasingstepsize} in a more general form than we need in this section. By \Cref{proxobjectivecombBG nonsmooth}, its hypothesis holds with $\rho_t \equiv 0$ under \eqref{eq: BG}. Under \eqref{eq: cw_cond} in \Cref{sec: componentwise_analysis}, however, the one-step inequality contains the extra term $\rho_t[h(\vx_t)-h^*]$, which we handle by choosing the weights of the convex combinations more carefully. This generalizes \citet[Lemma A.6]{vaidyan2026convergencerate} from constant step sizes and constant $\rho$ to non-increasing step sizes and time-varying $\rho_t$.

\begin{lemma}[Proximal last-iterate reduction]\label{theorem: proxsgd_decreasingstepsize}
    Let \Cref{assumption: prox} hold with positive non-increasing $(\tau_t)$, and define $S_t := \sum_{s=t}^T \tau_s$. Suppose there exist $0 \leq \rho_t \leq 1$ and $v_t \geq 0$ such that for all $0 \leq t \leq T$ and every $\mathcal{F}_t$-measurable $\vz_t$, we have
    \[
        \E[h(\vx_{t+1}) - h(\vz_t) - \rho_t h(\vx_t) + \rho_t h^*] \leq \frac{1}{2\tau_t} \E[\|\vx_t - \vz_t\|^2 - \|\vx_{t+1} - \vz_t\|^2] + v_t.
    \]
    Then, it holds that
    \[
        \E[h(\vx_{T+1}) - h^*]
        \leq \frac{1}{2} \exp\left( \sum_{t=1}^T \frac{\tau_t \rho_t}{S_{t}} \right) \left( \frac{ \|\vx_0 - \vx^*\|^2}{S_0} + 2 \sum_{t=0}^T v_t \frac{\tau_t}{S_{t}} + \frac{2 \rho_0\tau_0 [h(\vx_0) - h^*]}{S_0} \right).
    \]
\end{lemma}
The proof of \Cref{theorem: proxsgd_decreasingstepsize} is given in \Cref{app: proxsgd_decreasingstepsize_proof}.

As just discussed, the hypothesis of \Cref{theorem: proxsgd_decreasingstepsize} holds with $\rho_t=0$ due to \Cref{proxobjectivecombBG nonsmooth}, so the result of \Cref{theorem: proxsgd_decreasingstepsize} applies. We next instantiate this result with three standard step size choices.
\begin{theorem}[Proximal SGD under \eqref{eq: BG}]\label{theorem: proxsgd_bg_nonsmooth}
    Let \Cref{assumption: prox,assumption: oracle BG} hold. Then the following guarantees hold for the last iterate of \eqref{eq: prox_sgd_update}:
\begin{enumerate}[label=(\roman*)]
    \item \textbf{\textup{(Step size for $\log$-free rate)}} If $\tau_t = \tau\frac{T-t+1}{(T+1)^{3/2}}$ with $\tau = \frac{1}{4B}$, then for any $T \geq 1$, we have
    \[
    \E[h(\vx_{T+1})-h^*] = \mathcal{O}\left( \frac{B\|\vx_0-\vx^*\|^2 + G^2/B}{\sqrt{T+1}} \right).
    \]
    \item \textbf{\textup{(Anytime step size I)}} If $\tau_t = \frac{\tau}{\ln(t+2)\sqrt{t+2}}$ with $\tau = \frac{1}{4\sqrt{3}B}$, then for any $T\geq 2$, we have
    \[
    \E[h(\vx_{T+1})-h^*] = \mathcal{O}\left( \frac{B\|\vx_0-\vx^*\|^2 + G^2/B}{\sqrt{T+1}}\ln (T+2) \right).
    \]
    \item \textbf{\textup{(Anytime step size II)}} If $\tau_t = \frac{\tau}{\sqrt{t+1}}$ with $\tau < \frac{1}{4\sqrt{2}B}$, then for any $T \geq 1$, we have
    \[
        \E[h(\vx_{T+1})-h^*] = \mathcal{O}\left(\frac{\|\vx_0-\vx^*\|^2}{\tau\sqrt{T}}+\frac{\ln^2 (T+1)}{T^{1/2 - 16\tau^2 B^2}}\right).
    \]
\end{enumerate}
\end{theorem}
The proof of \Cref{theorem: proxsgd_bg_nonsmooth} is given in \Cref{app: BG_thm_nonsmooth}. The bounds with explicit constants for (i) and (ii) are given in \eqref{eq: proxSGD_nonsmooth_explicit_i} and \eqref{eq: proxSGD_nonsmooth_explicit_ii}, respectively. 
\begin{remark}\label{remark: projectedsgd_bg}
The result of \Cref{theorem: proxsgd_bg_nonsmooth} extends trivially to projected SGD for solving a constrained problem, under \eqref{eq: BG}, by setting $g=\delta_{\mathcal{C}}$.
\end{remark}
\begin{remark}\label{rmk: deterministic subgradient descent}
Specializing \Cref{remark: projectedsgd_bg} to the deterministic oracle $\widetilde{\nabla} f(\vx_t) := \vs_t \in \partial f(\vx_t)$ gives the (projected) subgradient method. In this case, \eqref{eq: BG} reduces to
\begin{equation}
\label{eq:det-BG}
    \| \vs \|^2 \le B^2 \| \vx - \vx_0\|^2 + G^2
    \qquad \text{for all } \vx \in \mathcal{C} \text{ and all }
    \vs \in \partial f(\vx),
\end{equation}
and \Cref{theorem: proxsgd_bg_nonsmooth} gives the same rates without expectations. This yields last-iterate guarantees for convex functions that are not necessarily smooth or Lipschitz; an example is given in \Cref{subsec: nonsmooth_example}.
\end{remark}

\paragraph{Rates.}
We now discuss the rates in \Cref{theorem: proxsgd_bg_nonsmooth}.
The horizon-dependent step size in (i) gives the rate $\mathcal{O}(1/\sqrt{T})$ with no logarithmic factors, which is the optimal order in $T$, but it requires knowing $T$ in advance.
The step size in (ii) does not depend on $T$, but this comes at the cost of an additional $\ln(T+2)$ factor.
The standard step size in (iii) is also independent of $T$, but it is not square summable. As discussed after \Cref{prop: BG boundedness}, the exponential factor therein grows as $T^{16\tau^2B^2}$ with this step size, which gives the slower rate.
The condition $\tau < 1/(4\sqrt{2}B)$ ensures $16\tau^2B^2 < 1/2$, so this upper bound still goes to zero.
Taking $\tau$ smaller brings the exponent closer to $1/2$, but increases the first term in the bound.

\paragraph{Comparison with prior work.}
Our results extend the works of \citet{garrigos2025last} and \citet{vaidyan2026convergencerate}, who required \eqref{eq: cw_cond}, which is stronger than our assumption in this section, and were limited to horizon-dependent step sizes. 
The other relevant work by \citet{liu2024revisitinglastiterateconvergencestochastic} handles variable step sizes but requires bounded variance and smoothness.
Finally, \Cref{theorem: proxsgd_bg_nonsmooth} (i), which uses the step size scheme of \citet{liu2024revisitinglastiterateconvergencestochastic}, removes the extra log factors that appear in the rates of \citet{garrigos2025last} and \citet{vaidyan2026convergencerate}. This answers an open question posed by \cite[Remark 3.6, Section 5]{garrigos2025last}.

\section{Applications of the B-G Assumption}\label{sec: examples}
In this section, we present examples to argue that the BG condition \eqref{eq: BG} applies to a broader class of problems than the component-wise convexity and smoothness \eqref{eq: cw_cond}. 
\subsection{Minimizing the Maximum of Quadratics}\label{subsec: nonsmooth_example} Consider the following nonsmooth problem given as
\begin{equation*}
    \min_{\vx \in \R^d} \phi(\vx) := \min_{\vx \in \R^d} \Big\{\max_{i = 1,\dots, n} \phi_i(\vx)\Big\},
\end{equation*}
where each $\phi_i$ is convex, differentiable and satisfies the BG condition \eqref{eq: BG}. That is, for all $\vx \in \R^d$, we assume
\begin{equation}\label{eq: BGMaxQuad}
    \|\grad \phi_i(\vx)\|^2 \leq B^2 \|\vx - \vx_0\|^2 + G^2.
\end{equation}
A concrete instance of this problem is the maximum of finitely many convex quadratics \citep[Section~5.3]{lu2019relative}:
\begin{equation*}
    \phi_i(\vx) = \frac{1}{2} \vx^\top A_i \vx + \vb_i^\top \vx + c_i,
\end{equation*} 
where $A_i$ is positive semi-definite. In this case, $\grad \phi_i(\vx) = A_i \vx + \vb_i = A_i(\vx - \vx_0) + (A_i \vx_0 + \vb_i)$, so each $\phi_i$ satisfies the BG condition since
\[
\begin{aligned}
\|\nabla \phi_i(\vx)\|^2 
&\leq 2 \max_{i=1,\dots,n} \|A_i\|^2 \|\vx - \vx_0\|^2 + 2 \max_{i=1,\dots,n} \|A_i \vx_0 + \vb_i\|^2.
\end{aligned}
\]
Notice that this problem is generally nonsmooth and not globally Lipschitz. Hence, most existing results for subgradient-based algorithms in the literature, such as the works of \citet{harvey2019tight,jain2019making,zamani2025exact}, do not apply. \citet{lu2019relative} solves this problem in the deterministic setting using mirror descent with a problem-specific generating function so that the resulting Bregman divergence controls the growth of the subgradients. This approach yields an $\mathcal{O}(1/\sqrt{T})$ convergence rate for the averaged iterate. In contrast, our analysis under \eqref{eq: BG} applies to the classical subgradient method, with no modification across different problems. Additionally, our analysis yields last-iterate guarantees, which \citet{lu2019relative} does not provide.

The subdifferential of $\phi$ is the convex hull of active gradients $\nabla \phi_i(\vx)$ where $i$ is in the set of active indices $I(\vx) = \arg\max_{j=1,\dots, n} \phi_j(\vx)$. 
An unbiased stochastic subgradient in this case can be obtained by sampling $i$ uniformly from the set of active indices and taking $\sgrad \phi(\vx)  = \grad \phi_i(\vx)$ as the stochastic subgradient since
\[
\begin{aligned}
\E [\sgrad \phi(\vx)] &= \frac{1}{|I(\vx)|} \sum_{i \in I(\vx)} \grad \phi_i(\vx) \in \partial \phi(\vx). \\
\end{aligned}
\]
Additionally, notice that $\phi$ satisfies the BG condition since
\[
\begin{aligned}
\E \|\sgrad \phi(\vx) \|^2 &= \frac{1}{|I(\vx)|} \sum_{i \in I(\vx)} \|\grad \phi_i(\vx)\|^2 \leq B^2 \|\vx - \vx_0\|^2 + G^2,
\end{aligned}
\]
where we used \eqref{eq: BGMaxQuad}. This example shows that the BG condition can handle nonsmooth and non-Lipschitz objectives arising from pointwise maxima. In particular, it provides a setting in which last-iterate guarantees for the classical stochastic subgradient method can be obtained through our BG framework, while \eqref{eq: cw_cond} does not apply. In fact, we are not aware of a last-iterate rate even for the \emph{deterministic} subgradient method for this problem.

\subsection{Black-Box Variational Inference}\label{subsec: bbvi}
In this section, we highlight another application of our results for black-box variational inference. Variational inference is a powerful technique for approximating an intractable posterior distribution, $p$, by using a tractable distribution from a given family.
Throughout this section, we follow the setting of \citet{domke2023provable} and refer to this work for a more detailed discussion of the variational formulation and Gaussian parameterizations.

Let $\vx$ denote the observed data and $\vz\in\R^d$ the latent variables, and let $p(\vx,\vz)$ be their joint density. The goal is to approximate the posterior $p(\vz\mid\vx)$ by a member of the Gaussian variational family
\[\begin{aligned}
    q_{\vw}(\vz)=\mathcal N(\vz\mid \vm,CC^\top),
    \qquad \vw=(\vm,C),
\end{aligned}
\]
where $\vm\in\R^d$ is the mean and $CC^\top$ is the covariance. This is done by minimizing the negative evidence lower bound (ELBO)
\begin{equation}\label{eq: ELBO}
    h(\vw)=f(\vw)+g(\vw),
\end{equation}
where we introduce the functions $f$ and $g$ next.

By using the reparameterization of $q_{\vw}(\vz)$ as $T_{\vw}(\vu)=\vm+C\vu$ with $\vu\sim\mathcal N(0,I)$, the energy term $f$ is given by
\begin{equation}\label{eq: energy_term}
    f(\vw) = \E_{\vu\sim\mathcal N(0,I)} \left[-\ln p(\vx,T_{\vw}(\vu))\right],
\end{equation}
and $g$ is the negative entropy term given by
\[\begin{aligned}
    g(\vw) = -\ln\det C-\frac d2\ln(2\pi e).
\end{aligned}\]
Our goal is to minimize $h$ over the parameters $\vw$ of the variational family. We assume that
\[\begin{aligned}
    \psi(\vz):=-\ln p(\vx,\vz)
\end{aligned}
\]
is convex and $L$-smooth and that the variational problem admits a solution $\vw^*$. Under these assumptions, $f$ is convex and $L$-smooth (see \citet[Lemmas 18 and 19]{domke2023provable}).

As in \citet[Section~3.1]{domke2023provable}, we will apply proximal SGD to minimize \eqref{eq: ELBO} over a domain we introduce next.
In particular, we restrict $C$ defining the family to be lower triangular with a positive diagonal, which gives us the optimization domain
\begin{equation*}
    \mathcal{W}_{\mathbb{T}} := \left\{ (\vm,C)\in\mathbb{R}^d\times\mathbb{T}^d: C_{ii}>0,\ i=1,\ldots,d \right\},
\end{equation*}
where $\mathbb{T}^d$ is the set of $d \times d$ lower-triangular matrices. Then, the optimization problem becomes
\[
    \min_{\vw\in\mathcal{W}_{\mathbb{T}}} \Big\{ f(\vw)+g(\vw):=\E_{\vu\sim\mathcal N(0,I)} \left[-\ln p(\vx,T_{\vw}(\vu))\right] -\ln\det C-\frac d2\ln(2\pi e)\Big\}.
\]
In this case, the unbiased estimator of the energy term is given by
\[\begin{aligned}
    \sgrad f(\vw;\vu) := -\grad_{\vw}\ln p(\vx,T_{\vw}(\vu)), \qquad \vu\sim\mathcal N(0,I),
\end{aligned}\]
and hence $\mathbb{E}_{\vu\sim \mathcal{N}(0, I)}[\sgrad f(\vw;\vu)] = \nabla f(\vw)$. By \citet[Theorem 2]{domke2023provable}, we have that
\begin{equation}\label{eq: energy_BG}
    \E\|\sgrad f(\vw;\vu)\|^2 \leq \beta_{\mathrm{energy}}^2\|\vw-\vw^*\|^2 + \gamma_{\mathrm{energy}}^2,
\end{equation}
for uniform constants  $\beta_{\mathrm{energy}}$, $\gamma_{\mathrm{energy}}$, and any $\vw^*$ in the solution set.

For $C\in\mathbb{T}^d$, \citet[Theorem 13]{domke2020provable} (see also \citealp[Lemma 35]{domke2023provable}) shows that the proximal operator for $g+\delta_{\mathcal{W}_{\mathbb{T}}}$ (needed to fit the problem to our composite template) is computed as

    \begin{equation*}
        \prox_{\tau_t (g+\delta_{\mathcal{W}_{\mathbb{T}}})}(\vw) = (\vm,C+\Delta C); \qquad \Delta C_{ii} = \frac{1}{2} \left(\sqrt{C_{ii}^2+4\tau_t} - C_{ii} \right).
    \end{equation*}

\begin{remark}\label{rmk: bbvi ratesa}
Applying \Cref{fact: bg centering} to \eqref{eq: energy_BG} with $\va = \vw^*$ and $\vb = \vw_0$ shows that $\sgrad f$ satisfies \eqref{eq: BG}. Hence, \Cref{theorem: proxsgd_bg_nonsmooth} applies, giving an anytime last-iterate rate of $\mathcal{O}(\ln T/\sqrt{T})$ for proximal SGD. Unlike \citet{domke2023provable,kim2025nearly}, we do not require strong convexity to get an anytime rate. When the horizon $T$ is known, the rate $\mathcal{O}(1/\sqrt{T})$ follows from \Cref{theorem: proxsgd_bg_nonsmooth}.
\end{remark}
\begin{remark}\label{rmk: bbvi ratesb}
Our results also cover the constrained formulation of the BBVI problem, which was described in \citet[Section~3.2]{domke2023provable}. In this setting, the stochastic gradient satisfies \eqref{eq: BG} \citep[Theorem~3]{domke2023provable}, so \Cref{remark: projectedsgd_bg} gives the same last-iterate rates for projected SGD.
\end{remark}

The key difference between the guarantees provided by \Cref{rmk: bbvi ratesa} and \Cref{rmk: bbvi ratesb}, and those of \citet{domke2023provable} is that our results provide last-iterate convergence guarantees in the convex setting, whereas \citet{domke2023provable} establish convergence guarantees for weighted averages of the iterates. In addition, our results provide anytime last-iterate guarantees using the step size scheme $\tau_t = \tau/(\ln(t+2)\sqrt{t+2})$,
whereas the parameters in \citet{domke2023provable} depend on a fixed horizon $T$. As stated in \citet[Section~2.4]{kim2025nearly}, it was an open question to derive anytime bounds for this setting, which our results in \Cref{theorem: proxsgd_bg_nonsmooth} and \Cref{remark: projectedsgd_bg} address. 

\section{Convergence Analysis under CW}\label{sec: componentwise_analysis}
In this section, we analyze the last-iterate convergence of projected and proximal SGD under the assumptions of component-wise convexity and smoothness in \eqref{eq: cw_cond}. 
Last-iterate guarantees in this setting were first obtained for SGD by \citet{garrigos2025last} and \citet{attia2025fast}. This analysis was then extended to proximal SGD by \citet{vaidyan2026convergencerate}. All of these analyses use the variance transfer bound \eqref{eq: ab_type} with a reference point $\vz_t$ in the style of \citet{zamani2025exact}. This gives the one-step inequality
\begin{equation*}
\E\left[h(\vx_{t+1}) - h(\vz_t)\right] \leq \frac{1}{2\tau_t}\E\left[\|\vx_t-\vz_t\|^2 - \|\vx_{t+1}-\vz_t\|^2\right] + 4\tau_t L\E\left[h(\vx_t)-h^*\right] + 2\tau_t\sigma_*^2.
\end{equation*}
In these works, $\tau_t = \tau$ is fixed. Compared with the bounded second moment case \citep[Eq. 6]{liu2024revisitinglastiterateconvergencestochastic}, there is an extra term $\rho\E[h(\vx_t)-h^*]$ with $\rho = 4\tau L$ on the right-hand side, and the weighted sum used by \citet{zamani2025exact} no longer cancels exactly. \citet{vaidyan2026convergencerate} handle the extra term with a constant step size $\tau_t = \tau$. In this case, $\rho$ does not depend on $t$, and the weights of the convex combinations defining $\vz_t$ can be adjusted to absorb $\rho\E[h(\vx_t)-h^*]$ \citep[see][Lemma A.6]{vaidyan2026convergencerate}. However, the argument relies on both the step size and $\rho$ being constant. It therefore does not cover horizon-free decreasing step sizes such as $\tau_t = \tau/\sqrt{t+1}$, which give anytime guarantees. It also does not cover the time-varying schedule of \citet{liu2024revisitinglastiterateconvergencestochastic}, which removes the logarithmic factor in the bounded variance setting.

Throughout this section, we require \Cref{assumption: prox,assumption: oracle CW}.
As in \Cref{sec: bg_analysis}, we follow the same roadmap, beginning with an inequality that controls the second moment of the stochastic gradients; see \citet[Lemma A.3]{khaled2023unified}.
\begin{lemma}[Variance transfer]
    Let \Cref{assumption: prox,assumption: oracle CW} hold. Then for every $\gamma > 0$ and $\vx\in\R^d$, we have
    \begin{equation}\label{variancetransfer}
        \E\|\grad f_i(\vx)\|^2 \leq 2L (1+ \gamma) [h(\vx) - h^*] + \left(1 + \frac{1}{\gamma}\right) \sigma_*^2.
    \end{equation}
\end{lemma}
Since every $f_i$ is $L$-smooth under \Cref{assumption: oracle CW}, it follows that $f$ is $L$-smooth. This allows us to use the following classical one-step inequality for proximal SGD.
\begin{lemma}
\label{lemma: onestepid_proximal}
    Let \Cref{assumption: prox,assumption: oracle CW} hold, and set $\tau_t \leq \frac{1}{2L}$.
    For all $t\geq 0$ and any $\mathcal{F}_t$-measurable $\vz_t$, we have
    \[\begin{aligned}
    2 \tau_t \E [h(\vx_{t+1}) - h(\vz_t)] \leq \E\|\vx_t - \vz_t\|^2 - \E \|\vx_{t+1} - \vz_t\|^2 + 2\tau_t^2 \E \|\sgrad f(\vx_t)\|^2.
    \end{aligned}\]
\end{lemma}
The proof of \Cref{lemma: onestepid_proximal} is given in \Cref{app: onestepid_proximal_proof}.
We now analyze one step of proximal SGD under \eqref{eq: cw_cond}.
\begin{lemma}[One-step of proximal SGD under \eqref{eq: cw_cond}]
\label{proxobjectivecomb_decreasingstepsize}
Let \Cref{assumption: prox,assumption: oracle CW} hold, and set $\tau_t \leq \frac{1}{4L}$ for all $t$.
    Then for all $t\geq 0$ and every $\mathcal{F}_t$-measurable $\vz_t$, \eqref{eq: prox_sgd_update} satisfies
    \begin{equation*}
        \begin{aligned}
            \E\left[ h(\vx_{t+1}) - h(\vz_t) - \rho_t h(\vx_t) + \rho_t h^* \right] &\leq
            \frac{1}{2\tau_t} \E\left[ \|\vx_t - \vz_t\|^2 - \|\vx_{t+1} - \vz_t\|^2 \right] + v_t,
        \end{aligned}
    \end{equation*}
    where $\rho_t = 4 \tau_t L $ and $v_t = 2\sigma_*^2 \tau_t$.
\end{lemma}
The proof of \Cref{proxobjectivecomb_decreasingstepsize} is given in \Cref{app: componentwise_lem1}. Using the above lemma with the proximal reduction in \Cref{theorem: proxsgd_decreasingstepsize}, we can obtain last-iterate guarantees for proximal SGD with decreasing step sizes. 
\begin{theorem}[Proximal SGD under \eqref{eq: cw_cond}]\label{theorem: convergenceratesforproximalsgd_componentwise}
Let \Cref{assumption: prox,assumption: oracle CW} hold, and set $\tau = \frac{1}{4L}$. Then the following guarantees hold for the last iterate of \eqref{eq: prox_sgd_update}:
\begin{enumerate}[label=(\roman*)]
    \item \textbf{\textup{(Step size for $\log$-free rate)}} If $\tau_t = \tau \frac{T-t+1}{(T+1)^{3/2}}$, then for any $T \geq 1$, we have
    \[
    \E[h(\vx_{T+1})-h^*] = \mathcal{O}\left( \frac{1}{\sqrt{T+1}}\left[ L\|\vx_0-\vx^*\|^2 + \frac{\sigma_*^2}{L} + \frac{h(\vx_0)-h^*}{T+1} \right] \right).
    \]
    \item \textbf{\textup{(Anytime step size)}} If $\tau_t = \frac{\tau}{\sqrt{t+1}}$, then for any $T \geq 1$, we have
    \[
    \E[h(\vx_{T+1})-h^*] = \mathcal{O}\left( \frac{1}{\sqrt{T+1}}\left[ L\|\vx_0-\vx^*\|^2 + \frac{\sigma_*^2\ln(T+1)}{L} + h(\vx_0)-h^* \right] \right).
    \]
\end{enumerate}
\end{theorem}
The proof of \Cref{theorem: convergenceratesforproximalsgd_componentwise} is given in \Cref{app: componentwise_thm1}. The bounds with explicit constants for (i) and (ii) are given in \eqref{eq: componentwise_explicit_i} and \eqref{eq: componentwise_explicit_ii}, respectively.

\begin{remark}\label{remark: projectedsgd_componentwise}
The result of \Cref{theorem: convergenceratesforproximalsgd_componentwise} extends trivially to projected SGD for solving a constrained problem, under \eqref{eq: cw_cond}, by setting $g=\delta_{\mathcal{C}}$.
\end{remark}

\paragraph{Rates.}
The horizon-dependent step size in \Cref{theorem: convergenceratesforproximalsgd_componentwise}(i) gives the rate $\mathcal{O}(1/\sqrt{T})$ with no logarithmic factors, which is the optimal order in $T$, but it requires knowing $T$ in advance.
The standard step size in (ii) does not depend on $T$, but this incurs an additional cost of a $\ln(T+1)$ factor, which multiplies the noise term $\sigma_*^2/L$.

\paragraph{Comparison with prior work.}
We extended the work of \citet{vaidyan2026convergencerate} by generalizing its analysis to a decreasing step size scheme that is not horizon-dependent.
Additionally, the step size in \Cref{theorem: convergenceratesforproximalsgd_componentwise}(i), which follows \citet{liu2024revisitinglastiterateconvergencestochastic}, removes the logarithmic factor, answering a question posed by \citet{garrigos2025last}.
This is a strict improvement over the results in \citet{vaidyan2026convergencerate}, which require a horizon-dependent step size scheme and obtain a rate of $\mathcal{O}(\ln(T)/\sqrt{T})$.

\begin{remark}\label{remark: oracle_comparison}
  While the BG condition \eqref{eq: BG} is the weaker assumption, this generality comes at the cost of step size flexibility. Under the BG condition, the standard anytime step size $\tau_t = \tau/\sqrt{t+1}$ only yields a suboptimal rate (see \Cref{theorem: proxsgd_bg_nonsmooth}(iii)), while under component-wise convexity and smoothness \eqref{eq: cw_cond} the same step size attains the optimal order up to a logarithmic factor (see \Cref{theorem: convergenceratesforproximalsgd_componentwise}(ii)).
\end{remark}

\section{Conclusion and Perspectives}\label{sec: conclusion}
In this paper, we studied last-iterate convergence rates of projected and proximal SGD under assumptions weaker than the classical uniformly bounded second moment condition \eqref{eq: bdd_grad_var}. We considered two different mechanisms for controlling the stochastic oracle. The first is based on the Blum-Gladyshev condition \eqref{eq: BG}, which controls the stochastic gradients through the distance to the initial iterate. The second is based on component-wise convexity and smoothness \eqref{eq: cw_cond}, which leads to ABC-type bounds \eqref{eq: ab_type} on the second moment of the stochastic gradients.
Our work opens up several directions for future research. For example, under the BG
condition, our guarantee for the standard anytime step size
$\tau_t=\tau/\sqrt{t+1}$ is suboptimal. It would be interesting to investigate whether this can be improved. Another direction is expanding our analysis to handle Bregman distances and establishing high probability guarantees on the last iterate.

\section*{AI Disclosure Statement}
After completing the proofs and the writing of the paper, GPT 6.0 (OpenAI) and Claude Opus 5.5 (Anthropic) were used for various copyediting tasks and for checking the proofs in this manuscript. For proof-checking, these tools were used to review derivations for errors such as incorrect indices, missing steps, or inconsistent constants. All suggested changes were reviewed and verified by the authors.

Our analysis of proximal SGD under \eqref{eq: BG} initially required $f$ to be smooth, and we had a separate section for projected SGD under \eqref{eq: BG} where $f$ was not necessarily smooth. 
When investigating whether smoothness could be removed for proximal SGD, we used these tools to explore whether such an extension was feasible, and their output suggested that it was. We then derived and wrote the proof of the nonsmooth result for proximal SGD, which allowed us to merge the separate sections for proximal and projected SGD. These tools were also used to assist in deriving the step size bounds in \Cref{proxboundsforstep4}. The resulting bounds and their proofs were verified by us.

\section*{Acknowledgments}
We acknowledge the support of the Natural Sciences and Engineering Research Council of Canada (NSERC), funding reference numbers RGPIN-2023-05250 (Friedlander) and RGPIN-2025-06634 (Alacaoglu).

\appendix
\crefalias{section}{appendix}\crefalias{subsection}{subappendix}
\section{Proofs for Section~\ref{sec: bg_analysis}}\label{sec: proofs_BG}
In this section, we provide the proofs of the various results in \Cref{sec: bg_analysis}. We start by establishing the boundedness of the iterates under \eqref{eq: BG}, which is a key step in our analysis.
\subsection{Proof of Proposition \ref{prop: BG boundedness}}\label{app: BG_prop1}
\begin{proof}
    We start from the hypothesis \eqref{eq: distance recursion} which states that
    \begin{equation}\label{eq: distance recursionproof}
    \E \|\vx_{t+1} - \vx^*\|^2 \leq \E\|\vx_t - \vx^*\|^2  + \tau_t^2 \E \|\vq_{t+1} - \sgrad f(\vx_t)\|^2.
    \end{equation}
    Using Young's inequality on the last term, we have
    \begin{equation}\label{eq: distance recursionproof2}
    \E \|\vq_{t+1} - \sgrad f(\vx_t)\|^2 \leq 2\E\|\vq_{t+1}\|^2 + 2\E\|\sgrad f(\vx_t)\|^2.
    \end{equation}
    We bound the second term on the right-hand side of \eqref{eq: distance recursionproof2} using the BG condition after recentering \eqref{eq: BG recentered}, which gives $\E\|\sgrad f(\vx_t)\|^2 \leq 2B^2\E\|\vx_t - \vx^*\|^2 + \widehat{G}^2$. For the first term on the right-hand side of \eqref{eq: distance recursionproof2}, we use Jensen's inequality and \eqref{eq: BG recentered} at $\vx_{t+1}$ to obtain
    \[\|\vq_{t+1}\|^2 = \|\E_{t+1}[\sgrad f(\vx_{t+1})]\|^2 \leq \E_{t+1}\|\sgrad f(\vx_{t+1})\|^2 \leq 2B^2\|\vx_{t+1} - \vx^*\|^2 + \widehat{G}^2.\]
    Plugging the above bounds into \eqref{eq: distance recursionproof2}, then back into \eqref{eq: distance recursionproof} and taking total expectation, we get
    \[\E \|\vx_{t+1} - \vx^*\|^2 \leq (1 + 4\tau_t^2 B^2) \E \|\vx_t - \vx^*\|^2 + 4\tau_t^2 B^2 \E \|\vx_{t+1} - \vx^*\|^2 + 4\tau_t^2 \widehat{G}^2.\]
    Rearranging and simplifying yields
    \[\begin{aligned}
    \E \|\vx_{t+1} - \vx^*\|^2 &\leq \frac{1 + 4\tau_t^2 B^2}{1 - 4\tau_t^2 B^2} \E \|\vx_t - \vx^*\|^2 + \frac{4\tau_t^2 \widehat{G}^2}{1 - 4\tau_t^2 B^2} \\
    &\leq (1 + 16\tau_t^2 B^2) \E \|\vx_t - \vx^*\|^2 + 8\tau_t^2 \widehat{G}^2,
    \end{aligned}\]
    where the second line used $\frac{1+a}{1-a} = 1 + \frac{2a}{1-a}$ and $\frac{1}{1-a} \leq 2$ where $a = 4\tau_t^2 B^2$ because $4\tau_t^2 B^2 \leq 1/2$.\footnote{Strictly speaking, we used here that $(\mathbb{E}\|\vx_t-\vx^*\|^2)$ is finite, which can be shown by induction.}
    Unrolling the recursion above, we get the bound
    \begin{equation}\label{eq: bgboundprop1}
        \E\|\vx_{t+1}-\vx^*\|^2 \leq \prod_{k=0}^t (1 + 16\tau_k^2 B^2) \|\vx_0 - \vx^*\|^2 + 8\sum_{k=0}^t \left( \prod_{j=k+1}^t (1 + 16\tau_j^2 B^2) \right) \tau_k^2 \widehat{G}^2.
    \end{equation}
    Now we bound the product on the right-hand side of this bound as
    \[\begin{aligned}
    \prod_{k=0}^t (1 + 16\tau_k^2 B^2) \leq \prod_{k=0}^t \exp(16\tau_k^2 B^2) = \exp\left(16\sum_{k=0}^t \tau_k^2 B^2\right) \leq \exp\left( 16 P_t B^2 \right),
    \end{aligned}\]
    where we used  $1+x \leq e^x$ and $\sum_{k=0}^t \tau_k^2 = P_t$. Substituting this back into \eqref{eq: bgboundprop1}, we get
    \begin{equation}\label{eq: bgboundprop2}
        \E\|\vx_{t+1}-\vx^*\|^2 \leq \exp\left( 16 P_t B^2 \right) \left[\|\vx_0 - \vx^*\|^2 +
        8 P_t \widehat{G}^2 \right].
    \end{equation}
    This is the first assertion. The second assertion follows immediately by plugging the bound on the iterates from \eqref{eq: bgboundprop2} into the BG condition \eqref{eq: BG recentered}.
\end{proof}
We next prove a one-step progress bound for proximal SGD under \eqref{eq: BG}.
\subsection{Proof of Lemma \ref{proxobjectivecombBG nonsmooth}}\label{app: BG_lem_nonsmooth}
\begin{proof}
    We apply \Cref{lemma: onestepid_proximal_nonsmooth} with $\vq_{t+1} := \E_{t+1}[\sgrad f(\vx_{t+1})] \in \partial f(\vx_{t+1})$ to obtain
    \begin{equation}\label{eq: simw3}
        2 \tau_t \E [h(\vx_{t+1}) - h(\vz_t)] \leq \E\|\vx_t - \vz_t\|^2 - \E \|\vx_{t+1} - \vz_t\|^2 + \tau_t^2 \E \|\vq_{t+1} - \sgrad f(\vx_t)\|^2.
    \end{equation}
    As in the proof of \Cref{prop: BG boundedness}, we use Young's inequality and Jensen's inequality to bound the last term as
    \begin{equation*}
        \E \|\vq_{t+1} - \sgrad f(\vx_t)\|^2 \leq 2\E\|\vq_{t+1}\|^2 + 2\E\|\sgrad f(\vx_t)\|^2 \leq 2\E\|\sgrad f(\vx_{t+1})\|^2 + 2\E\|\sgrad f(\vx_t)\|^2.
    \end{equation*}
    By \Cref{prop: BG boundedness}, we have that $\E\|\sgrad f(\vx_t)\|^2 \leq K_{t-1} \leq K_t$, since $P_t$ is non-decreasing, and also that $\E\|\sgrad f(\vx_{t+1})\|^2 \leq K_t$. Hence, the right-hand side is bounded by $4K_t$.
    Substituting this bound to \eqref{eq: simw3} and dividing by $2\tau_t$ gives us the desired result.
\end{proof}
We proceed to obtain the last-iterate convergence rate of proximal SGD under \eqref{eq: BG} with various step size schemes.
\subsection{Proof of Theorem \ref{theorem: proxsgd_bg_nonsmooth}}\label{app: BG_thm_nonsmooth}
\begin{proof}
    In all three step size schemes, $(\tau_t)$ is positive and non-increasing, and $\tau_t \leq \tau_0 \leq \frac{1}{2\sqrt{2}B}$. In particular, for (i), we have $\tau_0 = \frac{\tau}{\sqrt{T+1}} \leq \frac{1}{4B}$; for (ii), we have $\tau_0 = \frac{\tau}{\sqrt{2}\ln 2} \leq \frac{1}{2\sqrt{2}B}$ since $4\sqrt{3}\ln 2 \geq 2$; and for (iii), we have $\tau_0 = \tau < \frac{1}{4\sqrt{2}B}$. Hence, we can apply \Cref{proxobjectivecombBG nonsmooth} to get
    \begin{equation}\label{eq: tmp2 nonsmooth}
        \E [h(\vx_{t+1}) - h(\vz_t)] \leq \frac{1}{2\tau_t} \E[\|\vx_t - \vz_t\|^2 - \|\vx_{t+1} - \vz_t\|^2] + 2\tau_t K_t,
    \end{equation}
    where $K_t$ and $P_t$ are defined in \Cref{prop: BG boundedness}. Since $P_t$ is non-decreasing, we have that $P_t \leq P_T$, and hence, $K_t \leq K_T$. With this, we can upper bound the last term in \eqref{eq: tmp2 nonsmooth} by $2\tau_t K_T$.
    We can now apply \Cref{theorem: proxsgd_decreasingstepsize} with $\rho_t \equiv 0$ and $v_t = 2\tau_t K_T$ to obtain
\begin{equation}\label{eq: sieo4 nonsmooth}
    \E[h(\vx_{T+1})-h^*]
    \leq
    \frac{\|\vx_0-\vx^*\|^2}{2S_0}
    +
    2K_T\sum_{t=0}^T\frac{\tau_t^2}{S_t}.
\end{equation}
For the first two step size schemes, the corresponding bounds from
\Cref{proxboundsforstep2,proxboundsforstep4}, respectively, give
\begin{equation}\label{eq: proxSGD_nonsmooth_tab}
\begin{array}{c|c|c|c}
\text{Step size}
& S_0
& \displaystyle\sum_{t=0}^T\frac{\tau_t^2}{S_t}
& P_T \\[3ex]
\hline
\rule{0pt}{5ex}
\text{(i) } \tau_t = \tau\frac{T-t+1}{(T+1)^{3/2}}
& \displaystyle\geq\frac{\tau\sqrt{T+1}}{2}
& \displaystyle\leq\frac{2\tau}{\sqrt{T+1}}
& \leq\tau^2 \\[2ex]
\text{(ii) }  \tau_t = \frac{\tau}{\ln(t+2)\sqrt{t+2}}
& \displaystyle\geq
\frac{\tau\sqrt{T+1}}{\sqrt{2}\ln(T+2)}
& \displaystyle\leq
\frac{6\sqrt{2}\tau\ln(T+2)+33\tau}{\sqrt{T+1}}
& \leq3\tau^2 .
\end{array}
\end{equation}
Since $\tau = 1/(4B)$ in (i) and $\tau = 1/(4\sqrt{3}B)$ in (ii), the bounds on $P_T$ in \eqref{eq: proxSGD_nonsmooth_tab} give
$16P_TB^2\leq1$ for both cases. Using this bound in the definition of $K_T$ (in Prop. \ref{prop: BG boundedness}) and then expanding the definition of $\widehat G^2$ from \eqref{eq: BG recentered}, we get
\begin{equation}\label{eq: proxSGD_nonsmooth_KT_bound}
K_T \leq 2(2e+1)B^2 \|\vx_0 - \vx^*\|^2 + (1+e) G^2 =: K.
\end{equation}
Substituting the bounds from \eqref{eq: proxSGD_nonsmooth_tab} and \eqref{eq: proxSGD_nonsmooth_KT_bound} into \eqref{eq: sieo4 nonsmooth}, we obtain
\begin{align}
    \text{(i)}\quad \E[h(\vx_{T+1})-h^*]
    &\leq \frac{1}{\sqrt{T+1}}\left(\frac{\|\vx_0-\vx^*\|^2}{\tau} + 4\tau K\right), \label{eq: proxSGD_nonsmooth_explicit_i}\\
    \text{(ii)}\quad \E[h(\vx_{T+1})-h^*]
    &\leq \frac{1}{\sqrt{T+1}}\left[\frac{\ln(T+2)}{\sqrt{2}\tau}\|\vx_0-\vx^*\|^2 + 2\tau K\left(6\sqrt{2}\ln(T+2)+33\right)\right]. \label{eq: proxSGD_nonsmooth_explicit_ii}
\end{align}
In both cases, $\tau = \Theta(1/B)$. Substituting this and the definition of $K$ in \eqref{eq: proxSGD_nonsmooth_KT_bound} into \eqref{eq: proxSGD_nonsmooth_explicit_i} and \eqref{eq: proxSGD_nonsmooth_explicit_ii}, and suppressing constants, yields the first two claims of the theorem.

For the final step size scheme, \Cref{proxboundsforstep1} gives us the bounds
\[
S_0\geq2\tau(\sqrt{T+2}-1),
\qquad
\sum_{t=0}^T\frac{\tau_t^2}{S_t}
\leq\frac{2\tau(1+\ln(T+1))}{\sqrt{T+2}},
\qquad
P_T\leq\tau^2(1+\ln(T+1)).
\]
Substituting the above bounds into \eqref{eq: sieo4 nonsmooth}, we get
\begin{equation}\label{eq: proxSGD_BG_step3 nonsmooth}
\E[h(\vx_{T+1})-h^*]
\leq
\frac{\|\vx_0-\vx^*\|^2}
{4\tau(\sqrt{T+2}-1)}
+
K_T\frac{4\tau(1+\ln(T+1))}{\sqrt{T+2}}.
\end{equation}
Using the bound on $P_T$, we estimate $K_T$ (see Prop. \ref{prop: BG boundedness} for its definition) as follows:
\[
\begin{aligned}
K_T
&\leq
2B^2e^{16\tau^2B^2}(T+1)^{16\tau^2B^2}
\left(
\|\vx_0-\vx^*\|^2
+
8\tau^2(1+\ln(T+1))\widehat G^2
\right)
+\widehat G^2 \\
&=
\mathcal{O}\!\left(
T^{16\tau^2B^2}\ln (T+1)
\right).
\end{aligned}
\]
Substituting the bound on $K_T$ into \eqref{eq: proxSGD_BG_step3 nonsmooth} and suppressing constants give
\[
\E[h(\vx_{T+1})-h^*]
=
\mathcal{O}\!\left(
    \frac{\|\vx_0-\vx^*\|^2}{\tau\sqrt{T}}+
T^{-1/2+16\tau^2B^2}(\ln (T+1))^2
\right).
\]
In particular, for $\tau < 1/(4\sqrt{2}B)$, we have $-1/2 + 16\tau^2B^2 < 0$, and hence $\E[h(\vx_{T+1})-h^*] \to 0$ as $T \to \infty$. This concludes the proof of the theorem.
\end{proof}
\section{Proofs for Section~\ref{sec: componentwise_analysis}}\label{sec: proofs_componentwise}
In this section, we provide the proofs of the various results in \Cref{sec: componentwise_analysis}. We start by establishing a one-step progress bound for proximal SGD under \eqref{eq: cw_cond}.
\subsection{Proof of Lemma \ref{proxobjectivecomb_decreasingstepsize}}\label{app: componentwise_lem1}
\begin{proof}
    We start from the result of \Cref{lemma: onestepid_proximal} using $\vz = \vz_t$ to obtain,
    \[\begin{aligned}
    2 \tau_t \E [h(\vx_{t+1}) - h(\vz_t)] \leq \E\|\vx_t - \vz_t\|^2 - \E \|\vx_{t+1} - \vz_t\|^2 + 2\tau_t^2 \E \|\sgrad f(\vx_t)\|^2.
    \end{aligned}\]
    Dividing throughout by $2\tau_t$ and using \eqref{variancetransfer} with $\gamma = 1$, we get
    \[
    \begin{aligned}
        \E [h(\vx_{t+1}) - h(\vz_t)] &\leq \frac{1}{2\tau_t} \E[\|\vx_t - \vz_t\|^2 - \|\vx_{t+1} - \vz_t\|^2] + 4\tau_t L  \E[h(\vx_t) - h^*] + 2 \tau_t  \sigma_*^2.
    \end{aligned}
    \]
    Rearranging terms, we obtain the desired result.
\end{proof}
We now proceed to show the last-iterate convergence rate of proximal SGD under \eqref{eq: cw_cond} with various step size schemes.
\subsection{Proof of Theorem \ref{theorem: convergenceratesforproximalsgd_componentwise}}\label{app: componentwise_thm1}
\begin{proof}
From \Cref{proxobjectivecomb_decreasingstepsize}, we have for all
$t=0,\ldots,T$ and $\vz_t\in\mathcal{F}_t$,
\[
\E\left[
h(\vx_{t+1})-h(\vz_t)-\rho_t h(\vx_t)+\rho_t h^*
\right]
\leq
\frac{1}{2\tau_t}
\E\left[
\|\vx_t-\vz_t\|^2-\|\vx_{t+1}-\vz_t\|^2
\right]
+v_t,
\]
where $\rho_t=4\tau_tL$ and $v_t=2\sigma_*^2\tau_t$.
Since $\tau_t\leq\tau$ for both step size schemes and
$\tau L\leq1/4$, we have $\rho_t\leq1$. Hence,
\Cref{theorem: proxsgd_decreasingstepsize} gives
\begin{equation}\label{eq: componentwise_prox_master}
\E[h(\vx_{T+1})-h^*]
\leq
\frac{1}{2}
\exp\left(
\sum_{t=1}^T\frac{\tau_t\rho_t}{S_t}
\right)
\left(
\frac{\|\vx_0-\vx^*\|^2}{S_0}
+
4\sigma_*^2
\sum_{t=0}^T\frac{\tau_t^2}{S_t}
+
\frac{2\rho_0\tau_0[h(\vx_0)-h^*]}{S_0}
\right).
\end{equation}
We first bound the exponential term, by the definition $\rho_t=4\tau_tL$, as
\[
\exp\left(
\sum_{t=1}^T\frac{\tau_t\rho_t}{S_t}
\right)
=
\exp\left(
4L\sum_{t=1}^T\frac{\tau_t^2}{S_t}
\right).
\]
From
\Cref{proxboundsforstep1,proxboundsforstep2},
we obtain the following bounds for the two step size schemes:
\[
\begin{array}{c|c|c|c}
\text{Step size}
& S_0
& \displaystyle\sum_{t=0}^T\frac{\tau_t^2}{S_t}
& \displaystyle
\exp\left(
4L\sum_{t=1}^T\frac{\tau_t^2}{S_t}
\right)
\\[3ex]
\hline
\text{(i) } \tau_t = \tau\frac{T-t+1}{(T+1)^{3/2}}
& \displaystyle\geq\frac{\tau\sqrt{T+1}}{2}
& \displaystyle\leq\frac{2\tau}{\sqrt{T+1}}
& \leq e^2
\\[2ex]
\text{(ii) } \tau_t = \frac{\tau}{\sqrt{t+1}}
& \displaystyle\geq2\tau(\sqrt{T+2}-1)
& \displaystyle\leq
\frac{2\tau(1+\ln(T+1))}{\sqrt{T+2}}
& \leq e^{5/2}
\end{array}
\]
where we bounded the last column using $\tau L\leq1/4$ and the corresponding estimates for $\sum_{t=0}^T\frac{\tau_t^2}{S_t}$ from the third column.
Substituting these bounds into \eqref{eq: componentwise_prox_master}, and additionally using
\[
\frac{1}{\sqrt{T+2}-1}
\leq
\frac{2}{\sqrt{T+1}},
\qquad
\frac{1}{\sqrt{T+2}}
\leq
\frac{1}{\sqrt{T+1}}
\]
for (ii), we obtain
\begin{align}
    \text{(i)}\quad \E[h(\vx_{T+1})-h^*]
    &\leq \frac{e^2}{\sqrt{T+1}}\left[\frac{\|\vx_0-\vx^*\|^2}{\tau} + 4\tau\sigma_*^2 + \frac{8L\tau\,[h(\vx_0)-h^*]}{T+1}\right], \label{eq: componentwise_explicit_i}\\
    \text{(ii)}\quad \E[h(\vx_{T+1})-h^*]
    &\leq \frac{e^{5/2}}{\sqrt{T+1}}\left[\frac{\|\vx_0-\vx^*\|^2}{2\tau} + 4\tau\sigma_*^2(1+\ln(T+1)) + 4L\tau\,[h(\vx_0)-h^*]\right]. \label{eq: componentwise_explicit_ii}
\end{align}
Substituting $\tau = 1/(4L)$ into \eqref{eq: componentwise_explicit_i} and \eqref{eq: componentwise_explicit_ii}, and suppressing constants, yield the claims of the theorem.
\end{proof}

\section{Technical Lemmas for Proximal SGD}\label{sec: proximal_technical_lemmas}
We provide the proofs of the technical lemmas used in our analysis. We first restate a classical consequence of the definition of the proximal operator that we use often.
Its proof follows from the definition of $\prox_{\tau g}$ and the convexity of $g$ (see, e.g., \citealp[Theorem 6.39]{beck2017first}).
For $\vu=\prox_{\tau_t g}(\vy)$, it holds that
\begin{equation}\label{proxprojectionlemma}
    \langle \vu-\vy,\vx-\vu\rangle \geq \tau_t\bigl(g(\vu)-g(\vx)\bigr), \qquad \forall\,\vx\in\mathbb{R}^d.
\end{equation}
To bound the iterates and the second moment of the stochastic subgradients, we start from the following one-step inequality for \eqref{eq: prox_sgd_update}.
\subsection{Proof of Lemma \ref{lemma: onestepid_proximal_nonsmooth}}\label{app: onestepid_proximal_nonsmooth_proof}
\begin{proof}
    Using the identity $\langle \va, \vb \rangle = -\frac{1}{2} \| \va\|^2 - \frac{1}{2} \| \vb\|^2+\frac{1}{2} \| \va+\vb\|^2$, we have
    \begin{equation}\label{eq:prox1 nonsmooth}
        \|\vx_{t+1} - \vz_t\|^2 - \|\vx_t - \vz_t\|^2 = \|\vx_{t+1} - \vx_t\|^2 + 2\langle \vx_{t+1} - \vx_t, \vx_t - \vz_t\rangle.
    \end{equation}
    Consider the last term in the above identity. We use the prox-inequality \eqref{proxprojectionlemma} with $\vy= \vx_t - \tau_t \sgrad f(\vx_t)$ and $\vx = \vz_t$ to obtain
    \[
    \begin{aligned}
        \langle \vx_{t+1} - \vx_t, \vx_t - \vz_t \rangle &= \langle \vx_{t+1} - \vx_t, \vx_t - \vx_{t+1}\rangle + \langle \vx_{t+1} - \vx_t, \vx_{t+1} - \vz_t\rangle \\
        &\leq - \|\vx_{t+1} - \vx_{t}\|^2 + \tau_t (g(\vz_t) - g(\vx_{t + 1})) + \tau_t \langle \sgrad f(\vx_t), \vz_t - \vx_{t + 1}\rangle.
    \end{aligned}
    \]
    Substituting this into \eqref{eq:prox1 nonsmooth} gives
    \begin{equation}
        \|\vx_{t+1} - \vz_t\|^2 - \|\vx_t - \vz_t\|^2 \leq - \|\vx_{t+1} - \vx_t\|^2 + 2\tau_t (g(\vz_t) - g(\vx_{t + 1})) + 2\tau_t \langle \sgrad f(\vx_t), \vz_t - \vx_{t + 1}\rangle. \label{eq: spgd_term2 nonsmooth}
    \end{equation}
We can rewrite the last term as
\begin{equation}\label{eq: spgd_34 nonsmooth}
    \begin{aligned}
        2\tau_t \langle \sgrad f(\vx_t), \vz_t - \vx_{t + 1}\rangle 
        &= 2\tau_t \langle \sgrad f(\vx_t), \vz_t - \vx_{t}\rangle  + 2\tau_t \langle \sgrad f(\vx_t), \vx_{t} - \vx_{t + 1}\rangle.
    \end{aligned}
    \end{equation}
    Plugging this back to \eqref{eq: spgd_term2 nonsmooth} and taking conditional expectation with respect to $\mathcal{F}_t$, we have
    \begin{equation*}
        \begin{aligned}
            \E_t\|\vx_{t+1} - \vz_t\|^2 - \|\vx_t - \vz_t\|^2 &\leq - \E_t\|\vx_{t+1} - \vx_t\|^2 + 2\tau_t \E_t[g(\vz_t) - g(\vx_{t + 1})] \\
            &\quad + 2\tau_t \langle \vs_t, \vz_t - \vx_t\rangle + 2\tau_t \E_t\langle \sgrad f(\vx_t), \vx_t - \vx_{t + 1}\rangle,
        \end{aligned}
    \end{equation*}
    where $\vs_t := \E_t[\sgrad f(\vx_t)] \in \partial f(\vx_t)$ by \Cref{assumption: prox}, and we used that $\vx_t$ and $\vz_t$ are $\mathcal{F}_t$-measurable.
    We estimate $\langle \vs_t, \vz_t - \vx_t\rangle \leq [f(\vz_t) - f(\vx_t)]$ using the convexity of $f$. Adding and subtracting $2\tau_t f(\vx_{t+1})$ and rearranging terms give us
        \begin{align}\label{eq: spgd_5}
            \E_t\|\vx_{t+1} - \vz_t\|^2 - \|\vx_t - \vz_t\|^2 &\leq 2\tau_t \E_t[h(\vz_t) - h(\vx_{t + 1})] - \E_t\|\vx_{t+1} - \vx_t\|^2 \notag\\
            &\quad  + 2\tau_t \E_t[f(\vx_{t+1}) - f(\vx_t)] + 2\tau_t \E_t\langle \sgrad f(\vx_t), \vx_t - \vx_{t + 1}\rangle.
        \end{align}
    We now estimate the last two terms in the above inequality. Let $\vq_{t+1} \in \partial f(\vx_{t+1})$. By convexity, we have $f(\vx_{t+1}) - f(\vx_t) \leq \langle \vq_{t+1}, \vx_{t+1} - \vx_t\rangle$.
    Using this to bound the second line in \eqref{eq: spgd_5} gives
    \begin{equation}\label{eq: spgd_6}
        2\tau_t [f(\vx_{t+1}) - f(\vx_t)] + 2\tau_t \langle \sgrad f(\vx_t), \vx_t - \vx_{t + 1}\rangle \leq 2\tau_t \langle \vq_{t+1} - \sgrad f(\vx_t), \vx_{t+1} - \vx_t\rangle.
    \end{equation}
    We estimate the right-hand side using Young's inequality as
    \begin{equation}\label{eq: spgd_7}
        2\tau_t \langle \vq_{t+1} - \sgrad f(\vx_t), \vx_{t+1} - \vx_t\rangle \leq \tau_t^2 \|\vq_{t+1} - \sgrad f(\vx_t)\|^2 +  \|\vx_{t+1} - \vx_t\|^2.
    \end{equation}
    We combine \eqref{eq: spgd_6} and \eqref{eq: spgd_7} (after taking conditional expectation) in \eqref{eq: spgd_5} and group the terms involving $\|\vx_{t+1} - \vx_t\|^2$ to obtain
    \begin{equation*}
        \begin{aligned}
            \E_t\|\vx_{t+1} - \vz_t\|^2 - \|\vx_t - \vz_t\|^2 &\leq 2\tau_t \E_t[h(\vz_t) - h(\vx_{t + 1})] + \tau_t^2 \E_t\|\vq_{t+1} - \sgrad f(\vx_t)\|^2.
        \end{aligned}
    \end{equation*}
    Taking total expectation gives us the assertion.
\end{proof}
We next prove the analogue of \Cref{lemma: onestepid_proximal_nonsmooth} under smoothness.
\subsection{Proof of Lemma \ref{lemma: onestepid_proximal}}\label{app: onestepid_proximal_proof}
\begin{proof}
    Using the identity $\langle \va, \vb \rangle = -\frac{1}{2} \| \va\|^2 - \frac{1}{2} \| \vb\|^2+\frac{1}{2} \| \va+\vb\|^2$, we have
    \begin{equation}\label{eq:prox1}
        \|\vx_{t+1} - \vz_t\|^2 - \|\vx_t - \vz_t\|^2 = \|\vx_{t+1} - \vx_t\|^2 + 2\langle \vx_{t+1} - \vx_t, \vx_t - \vz_t\rangle.
    \end{equation}
    Consider the last term in the above identity. We use the prox-inequality \eqref{proxprojectionlemma} with $\vy= \vx_t - \tau_t \sgrad f(\vx_t)$ and $\vx = \vz_t$ to obtain
    \[
    \begin{aligned}
        \langle \vx_{t+1} - \vx_t, \vx_t - \vz_t \rangle &= \langle \vx_{t+1} - \vx_t, \vx_t - \vx_{t+1}\rangle + \langle \vx_{t+1} - \vx_t, \vx_{t+1} - \vz_t\rangle \\
        &\leq - \|\vx_{t+1} - \vx_{t}\|^2 + \tau_t (g(\vz_t) - g(\vx_{t + 1})) + \tau_t \langle \sgrad f(\vx_t), \vz_t - \vx_{t + 1}\rangle.
    \end{aligned}
    \]
    Substituting this into \eqref{eq:prox1} gives
    \begin{equation}
        \|\vx_{t+1} - \vz_t\|^2 - \|\vx_t - \vz_t\|^2 \leq - \|\vx_{t+1} - \vx_t\|^2 + 2\tau_t (g(\vz_t) - g(\vx_{t + 1})) + 2\tau_t \langle \sgrad f(\vx_t), \vz_t - \vx_{t + 1}\rangle. \label{eq: spgd_term2}
    \end{equation}
We can rewrite the last term as
\begin{equation}\label{eq: spgd_34}
    \begin{aligned}
        2\tau_t \langle \sgrad f(\vx_t), \vz_t - \vx_{t + 1}\rangle 
        &= 2\tau_t \langle \sgrad f(\vx_t), \vz_t - \vx_{t}\rangle + 2\tau_t \langle \sgrad f(\vx_t) - \grad f(\vx_t), \vx_t - \vx_{t + 1}\rangle \\
        &\quad  + 2\tau_t \langle \grad f(\vx_t), \vx_t - \vx_{t + 1}\rangle.
    \end{aligned}
    \end{equation}
    We estimate the second term in \eqref{eq: spgd_34} using Young's inequality as
    \begin{equation*}
        2 \tau_t \langle \sgrad f(\vx_t) - \grad f(\vx_t), \vx_t - \vx_{t+1} \rangle\leq \frac{\tau_t^2}{1-\tau_tL} \|\sgrad f(\vx_t) - \grad f(\vx_t)\|^2 + (1-\tau_t L) \| \vx_t-\vx_{t+1}\|^2.
    \end{equation*}
    For the last term in \eqref{eq: spgd_34}, we can use the smoothness of $f$ to derive
    \begin{equation*}
        2\tau_t \langle \grad f(\vx_t), \vx_t - \vx_{t + 1}\rangle \leq 2\tau_t \Big( f(\vx_t) - f(\vx_{t+1}) + \frac{L}{2} \|\vx_t - \vx_{t+1}\|^2\Big).
    \end{equation*}
    We plug the last two estimates into \eqref{eq: spgd_34}, which we then insert into \eqref{eq: spgd_term2} and rearrange to obtain
    \begin{equation*}
        \begin{aligned}
                \|\vx_{t+1} - \vz_t\|^2 - \|\vx_t - \vz_t\|^2  & \leq 2\tau_t \langle \sgrad f(\vx_t), \vz_t - \vx_{t}\rangle + 2\tau_t (g(\vz_t) - g(\vx_{t + 1})) \\
        & \quad + \frac{\tau_t^2}{1-\tau_t L} \|\sgrad f(\vx_t) - \grad f(\vx_t)\|^2 +2\tau_t [f(\vx_t) - f(\vx_{t+1})].
        \end{aligned}
    \end{equation*}
    Since $\tau_t \leq \frac1{2L}$, we have $\frac{\tau_t^2}{1-\tau_t L} \leq 2\tau_t^2$. Plugging in this estimate, and taking conditional expectation with respect to $\mathcal{F}_t$, we obtain,
    \[\begin{aligned}
    \E_t\|\vx_{t+1} - \vz_t\|^2 &\leq \|\vx_t - \vz_t\|^2 + 2\tau_t \langle \grad f(\vx_t), \vz_t - \vx_t \rangle + 2 \tau_t \E_t[g(\vz_t) - g(\vx_{t+1})] \\
    & \quad + 2\tau_t^2 \E_t[\|\sgrad f(\vx_t) - \grad f(\vx_t)\|^2] + 2\tau_t \E_t[f(\vx_t) - f(\vx_{t+1})].
    \end{aligned}\]
    Using  $\E_t[\|\sgrad f(\vx_t) - \grad f(\vx_t)\|^2] \leq \E_t \|\sgrad f(\vx_t)\|^2$, the convexity of $f$, and $h = f + g$, we obtain
    \[\begin{aligned}
        \E_t\|\vx_{t+1} - \vz_t\|^2 &\leq \|\vx_t - \vz_t\|^2 + 2\tau_t \E_t[h(\vz_t) - h(\vx_{t+1})] + 2\tau_t^2 \E_t[\|\sgrad f(\vx_t)\|^2]. \\
    \end{aligned}\]
    Taking total expectation and rearranging gives the desired result.    
\end{proof}
    
We next prove the last-iterate reduction for proximal SGD, which converts the one-step bound into a bound on the last iterate.
\subsection{Proof of Lemma \ref{theorem: proxsgd_decreasingstepsize}}\label{app: proxsgd_decreasingstepsize_proof}
\begin{proof}
    Let us define $\alpha_{-1} = \alpha_0 = 1$ and $\alpha_t = \frac{S_{t-1}}{S_t + \tau_t \rho_t}\alpha_{t-1}$ for $t = 1, \ldots, T$. Since $\rho_t \leq 1$ and $\tau_{t-1} \geq \tau_t$, we have $S_{t-1} = S_t + \tau_{t-1} \geq S_t + \tau_t\rho_t$. As a result, we have that the sequence $(\alpha_t)$ is non-decreasing. Let $\vz_{-1} = \vx^*$ and
    \[
        \vz_t = \left(1 - \frac{\alpha_{t-1}}{\alpha_t}\right)\vx_t + \frac{\alpha_{t-1}}{\alpha_t}\vz_{t-1}, \qquad t = 0, \ldots, T.
    \]
    Unrolling the definition of $\vz_t$ shows that it is a convex combination of $\vx^*, \vx_0, \ldots, \vx_t$, since the coefficients below are nonnegative and sum to $1$:
    \begin{equation}\label{eq: ztprox}
    \vz_t = \frac{1}{\alpha_t} \vx^* + \frac{1}{\alpha_t} \sum_{s=0}^t (\alpha_s - \alpha_{s-1}) \vx_s.
    \end{equation}
    Hence $\vz_t$ is $\mathcal{F}_t$-measurable, so the hypothesis of this lemma applies with this $\vz_t$. Multiplying both sides of the hypothesis by $2\alpha_t\tau_t$ gives
\begin{equation}\label{eq: proxcombdec131}
\begin{aligned}
    2 \alpha_t \tau_t \E [h(\vx_{t+1}) -  h(\vz_t) - \rho_t h(\vx_t) + \rho_t h^*] &\leq \alpha_t \E[\|\vx_t - \vz_t\|^2 - \|\vx_{t+1} - \vz_t\|^2] + 2 \alpha_t v_t \tau_t. 
\end{aligned}
\end{equation}
The definition of $\vz_t$ also gives
    \begin{equation*}
    \vx_t - \vz_t = \vx_t - \left( 1-\frac{\alpha_{t-1}}{\alpha_t} \right)\vx_t - \frac{\alpha_{t-1}}{\alpha_t}\vz_{t-1} = \frac{\alpha_{t-1}}{\alpha_t}\left(\vx_t-\vz_{t-1} \right).
    \end{equation*}
Using this and $\alpha_{t-1} \leq \alpha_t$, we estimate the first term on the right-hand side of \eqref{eq: proxcombdec131} as
    \begin{equation*}
    \alpha_t \E\|\vx_t-\vz_t\|^2 = \frac{\alpha_{t-1}^2}{\alpha_t} \E\| \vx_t-\vz_{t-1}\|^2 \leq \alpha_{t-1}\E\| \vx_t-\vz_{t-1}\|^2.
    \end{equation*}
Substituting this estimate into \eqref{eq: proxcombdec131} and dividing by $2$, we get
\[
\begin{aligned}
    &\alpha_t \tau_t \E [h(\vx_{t+1}) - h(\vz_t) - \rho_t h(\vx_t) + \rho_t h^*] \\
    &\qquad \leq \frac{\alpha_{t-1}}{2} \E\|\vx_t - \vz_{t-1}\|^2 - \frac{\alpha_t}{2} \E\|\vx_{t+1} - \vz_t\|^2 + \alpha_t v_t \tau_t.
\end{aligned}
\]
Summing from $t = 0$ to $T$ gives
\begin{equation}\label{proxcombdec114}
    \Gamma:= \sum_{t=0}^T \alpha_t \tau_t \E [h(\vx_{t+1}) -  h(\vz_t) - \rho_t h(\vx_t) + \rho_t h^*] 
    \leq \frac{1}{2} \|\vx_0 - \vx^*\|^2 + \sum_{t=0}^T \alpha_t v_t \tau_t,
\end{equation}
where we dropped the negative term on the right-hand side and used the fact $\vz_{-1} = \vx^*$ and $\alpha_{-1} = 1$.
    We now lower bound $\Gamma$. By the convexity of $h$ applied to \eqref{eq: ztprox}, we have
    \begin{align*}
    h(\vz_t) \leq \frac{1}{\alpha_t} h(\vx^*) + \frac{1}{\alpha_t}\sum_{s=0}^t (\alpha_s - \alpha_{s-1}) h(\vx_s),
    \end{align*}
    and consequently
\begin{equation}\label{proxcombdec115}
    \begin{aligned}
        \Gamma
        & \geq  \sum_{t=0}^{T} \tau_t \E\left[
        \alpha_t h(\vx_{t+1})
        - \sum_{s=0}^t (\alpha_s - \alpha_{s-1}) h(\vx_s)
        - \rho_t \alpha_t h(\vx_t)
        + \rho_t\alpha_t  h^* - h^*
        \right].
    \end{aligned}
    \end{equation}
    Now consider the second term inside the expectation. We rewrite it as follows:
    \[\begin{aligned}
        \sum_{s=0}^t (\alpha_s - \alpha_{s-1}) h(\vx_s) 
        &= \sum_{s=0}^t (\alpha_s - \alpha_{s-1})\left[h(\vx_s) - h^*\right] + (\alpha_t - 1)h^*,
    \end{aligned}
    \]
    where we used the fact $\sum_{s=0}^t (\alpha_s - \alpha_{s-1}) = \alpha_t - \alpha_{-1} = \alpha_t - 1$. Substituting this back into \eqref{proxcombdec115} and simplifying gives
    \begin{equation}
    \begin{aligned}
        \Gamma \geq
        \sum_{t=0}^{T} \tau_t \E\left[
        \alpha_t (h(\vx_{t+1}) - h^*) - \rho_t \alpha_t(h(\vx_t) - h^*) - \sum_{s=0}^t (\alpha_s - \alpha_{s-1})\left[h(\vx_s) - h^*\right]
        \right]. \label{proxcombdec117}
    \end{aligned}
\end{equation}
    We next proceed to show that this lower bound gives us an expression for the last iterate. We define 
    \begin{equation*}
    r_t = h(\vx_t) - h^*.
    \end{equation*}
     Then \eqref{proxcombdec117} becomes
\begin{equation}\label{eq: swrdec4}
    \begin{aligned}
        \Gamma \geq \E\left[
        \sum_{t=0}^{T} \tau_t \alpha_t r_{t+1}
        - \sum_{t=0}^{T} \rho_t \tau_t \alpha_t r_t
        - \sum_{t=0}^{T}\sum_{s=0}^{t} \tau_t (\alpha_s - \alpha_{s-1}) r_s
        \right].
    \end{aligned}
    \end{equation}
    Consider the first term on the right-hand side of \eqref{eq: swrdec4}. We can rewrite it as
    \[
    \begin{aligned}
        \sum_{t=0}^{T} \tau_t \alpha_t r_{t+1}
        &= \tau_T \alpha_T r_{T+1} + \sum_{t=1}^{T} \tau_{t-1} \alpha_{t-1} r_t.
    \end{aligned}
    \]
    The second term on the right-hand side of \eqref{eq: swrdec4} can be rewritten as
    \[ \sum_{t=0}^{T} \rho_t \tau_t \alpha_t r_t = \rho_0 \tau_0 r_0 + \sum_{t=1}^{T} \rho_t \tau_t \alpha_t r_t,\]
    where we used the fact that $\alpha_0 = 1$. And finally, for the last term on the right hand side of \eqref{eq: swrdec4}, we swap the order of summation and use the definition of $S_t$ to get
\[
\sum_{t=0}^{T} \sum_{s=0}^{t} \tau_t(\alpha_s-\alpha_{s-1})r_s = \sum_{t=0}^{T} (\alpha_t-\alpha_{t-1})r_t S_t = \sum_{t=1}^{T} (\alpha_t - \alpha_{t-1})r_t S_t,
\]
where we used the fact that $\alpha_{-1} = \alpha_0 = 1$ so the first term ($t=0$) in the sum is zero.

    Combining the last three display equations in \eqref{eq: swrdec4} and grouping terms involving $r_t$ gives
\begin{equation}\label{eq:sgf4}
\Gamma \geq \E \left[
    \tau_T \alpha_T r_{T+1} + \sum_{t=1}^{T} r_t \left(
        \tau_{t-1} \alpha_{t-1}
        - (\alpha_t - \alpha_{t-1}) S_t
        - \tau_t \rho_t \alpha_t
    \right) - \rho_0 \tau_0 r_0
\right].
\end{equation}
From the definition of $\alpha_t$ and $S_{t-1} = S_t + \tau_{t-1}$, we have
\[
\begin{aligned}
    \tau_{t-1} \alpha_{t-1} - (\alpha_t - \alpha_{t-1}) S_t - \tau_t \rho_t \alpha_t &= 0.
\end{aligned}
\]
Using this in \eqref{eq:sgf4}, and plugging the result to \eqref{proxcombdec114}, we obtain the following inequality after moving $\rho_0 \tau_0 r_0$ to the right-hand side and dividing by $\tau_T \alpha_T$:
\begin{equation}
\begin{aligned}
    \E[r_{T+1}] &\leq \frac{1}{2}\left(\frac{ \|\vx_0 - \vx^*\|^2}{\tau_T \alpha_T} + 2 \sum_{t=0}^T \frac{\alpha_t v_t \tau_t}{\tau_T \alpha_T} + \frac{2\rho_0 \tau_0 r_0}{\tau_T \alpha_T}\right).
\end{aligned}\label{eq: proxcombdec118}
\end{equation}
To estimate the factor, $\frac{1}{\tau_T \alpha_T}$, let 
\[ \Lambda_T := \prod_{t=1}^T \left(1+\frac{\tau_t \rho_t}{S_t}\right). \]
By the definition of $S_t$ in the statement of the lemma, we have that $S_T = \tau_T$. Using this and  the definition of $\alpha_t$ (along with $\alpha_0=1$) give
\[
\begin{aligned}
    \frac{1}{\tau_T \alpha_T} &= \frac{1}{S_T} \prod_{t=1}^T\frac{\alpha_{t-1}}{\alpha_t}= \frac{1}{S_T} \prod_{t=1}^T \frac{S_t + \tau_t \rho_t}{S_{t-1}} = \frac{1}{S_T} \prod_{t=1}^T  \frac{S_t}{S_{t-1}} \left(1 + \frac{\tau_t \rho_t}{S_{t}} \right) = \frac{1}{S_0} \Lambda_T.
\end{aligned}
\]
We next estimate the second term inside the parenthesis in \eqref{eq: proxcombdec118} similarly:
\begin{equation} \label{eq: proxcombdec211}
\begin{aligned}
    \frac{\alpha_t v_t \tau_t}{\tau_T \alpha_T} &= \frac{v_t \tau_t}{S_T} \prod_{i=t+1}^T \frac{S_i + \tau_i \rho_i}{S_{i-1}} = \frac{v_t \tau_t}{S_T} \prod_{i=t+1}^T  \frac{S_i}{S_{i-1}} \left(1 + \frac{\tau_i \rho_i}{S_{i}} \right) \leq \frac{v_t \tau_t}{S_{t}} \Lambda_T.
\end{aligned}
\end{equation}
Plugging these estimates into \eqref{eq: proxcombdec118}, we get
\begin{equation}\label{eq: proxcombdec119}
\begin{aligned}
    \E[r_{T+1}] &\leq \frac{1}{2} \Lambda_T \left( \frac{ \|\vx_0 - \vx^*\|^2}{S_0} + 2 \sum_{t=0}^T v_t \frac{\tau_t}{S_{t}} + \frac{2 \rho_0 \tau_0 r_0}{S_0} \right). \\
\end{aligned}
\end{equation}
Now using the fact that $1+x \leq e^x$ for all $x$, we can estimate $\Lambda_T$ as follows:
\[\begin{aligned}
    \Lambda_T &= \prod_{t=1}^T \left(1+\frac{\tau_t \rho_t}{S_t}\right) \leq \prod_{t=1}^T \exp\left(\frac{\tau_t \rho_t}{S_t}\right) =\exp\left(\sum_{t=1}^T \frac{\tau_t \rho_t}{S_t}\right).
\end{aligned}
\]
Substituting this back into \eqref{eq: proxcombdec119} yields the assertion.
\end{proof}

\subsection{Technical Lemmas for Bounding the Step Size Terms}\label{sec: stepsize_technical_lemmas}
In this section, we prove various bounds related to the step sizes. Throughout this section, we define $S_t := \sum_{s=t}^T \tau_s$.
\begin{lemma}\label{proxboundsforstep1}
Let us set $\tau_t=\frac{\tau}{\sqrt{t+1}}$ for some $\tau>0$.
Then, for any $T\geq1$, we have
{\small
\[
\begin{alignedat}{2}
\text{(i)}\;\; \sum_{t=0}^T \tau_t^2
&\leq \tau^2(1+\ln(T+1)),
\qquad &
\text{(ii)}\;\; S_0
&\geq 2\tau(\sqrt{T+2}-1),
\\
\text{(iii)}\;\; \sum_{t=0}^T \frac{\tau_t^2}{S_t}
&\leq \frac{2\tau(1+\ln(T+1))}{\sqrt{T+2}},
\qquad &
\text{(iv)}\;\; \sum_{t=0}^{T-1} \frac{\tau_t^2}{S_{t+1}}
&\leq \frac{2\tau(1+\ln T)}{\sqrt{T}}.
\end{alignedat}
\]
}
\end{lemma}
\begin{proof}
Assertion (i) follows trivially by summing the harmonic series.
Next, by the definition of $S_t$ in the statement of \Cref{theorem: proxsgd_decreasingstepsize}, we have
\[
\begin{aligned}
S_t = \sum_{s=t}^T \frac{\tau}{\sqrt{s+1}} \geq \int_t^{T+1} \frac{\tau}{\sqrt{x+1}} dx = 2\tau(\sqrt{T+2} - \sqrt{t+1}),
\end{aligned}
\]
which gives (ii) after setting $t=0$.

We now estimate (iii) by using the lower bound for $S_t$, just derived, and the definition of $\tau_t$ as
\begin{equation}\label{eq: soxu4}
\begin{aligned}
\sum_{t=0}^T \frac{\tau_t^2}{S_t}
\leq \tau\sum_{t=0}^T
\frac{1}{2(t+1)(\sqrt{T+2} - \sqrt{t+1})} 
\leq \tau\sqrt{T+2}\sum_{t=0}^T
\frac{1}{(t+1)(T+1-t)},
\end{aligned}
\end{equation}
since 
\begin{align*}
\frac{1}{\sqrt{T+2}-\sqrt{t+1}} 
= \frac{\sqrt{T+2}+\sqrt{t+1}}{T+1-t} 
\leq \frac{2\sqrt{T+2}}{T+1-t}.
\end{align*}
Now, on \eqref{eq: soxu4}, we use
\[ \frac{1}{(t+1)(T+1-t)} = \frac{1}{T+2} \left( \frac{1}{t+1} + \frac{1}{T+1-t} \right), \]
to obtain
\[ \begin{aligned} 
    \sum_{t=0}^T \frac{\tau_t^2}{S_t} &\leq \frac{\tau}{\sqrt{T+2}} \sum_{t=0}^T \left( \frac{1}{t+1} + \frac{1}{T+1-t} \right) = \frac{2\tau}{\sqrt{T+2}} \sum_{t=0}^T \frac{1}{t+1} \leq \frac{2\tau(1+\ln(T+1))}{\sqrt{T+2}}. 
    \end{aligned} 
\]
Finally, assertion (iv) follows from a similar estimation. Since
\[\begin{aligned}
\sum_{t=0}^{T-1} \frac{\tau_t^2}{S_{t+1}} 
\leq \frac{\tau}{2} \sum_{t=0}^{T-1} \frac{1}{(t+1)(\sqrt{T+2} - \sqrt{t+2})} 
\leq \tau\sqrt{T+2} \sum_{t=0}^{T-1} \frac{1}{(t+1)(T-t)},
\end{aligned}\]
we repeat the same argument as in the proof of the third assertion by replacing $T+1-t$ with $T-t$ to get
\begin{align*}
\sum_{t=0}^{T-1} \frac{\tau_t^2}{S_{t+1}} \leq \frac{2\tau\sqrt{T+2}(1 + \ln(T))}{T+1}.
\end{align*}
The claim then follows by $\frac{\sqrt{T+2}}{T+1} \leq \frac{1}{\sqrt{T}}$.
\end{proof}
\begin{lemma}\label{proxboundsforstep2}
Let us set $\tau_t=\tau\frac{T-t+1}{(T+1)^{3/2}}$ for some $\tau>0$. Then, for any $T\geq1$, we have
{\small
\[
\begin{alignedat}{2}
\text{(i)}\;\; \sum_{t=0}^T \tau_t^2
&\leq \tau^2,
\qquad &
\text{(ii)}\;\; S_0
&\geq \frac{\tau\sqrt{T+1}}{2},
\\
\text{(iii)}\;\; \sum_{t=0}^T \frac{\tau_t^2}{S_t}
&\leq \frac{2\tau}{\sqrt{T+1}},
\qquad &
\text{(iv)}\;\; \sum_{t=0}^{T-1} \frac{\tau_t^2}{S_{t+1}}
&\leq \frac{4\tau}{\sqrt{T+1}}.
\end{alignedat}
\]
}
\end{lemma}
\begin{proof}
By the definition of $\tau_t$, we have
\begin{equation}\label{eq: fol4}
    \sum_{s=0}^T \tau_s^2 = \frac{\tau^2}{(T+1)^3} \sum_{s=0}^T (T-s+1)^2.
    \end{equation} 
We estimate the sum as
\begin{align*}
\sum_{s=0}^T (T-s+1)^2 = \sum_{k=1}^{T+1} k^2 = \frac{(T+1)(T+2)(2T+3)}{6}.
\end{align*}
Plugging this into \eqref{eq: fol4} and using $(T+2)(2T+3)\leq 6(T+1)^2$ gives assertion (i).

For assertion (ii), we estimate $S_t$ (defined in the statement of \Cref{theorem: proxsgd_decreasingstepsize}) as
\begin{equation}\label{eq: stbd1} 
    S_t = \sum_{s=t}^T \tau_s = \frac{\tau}{(T+1)^{3/2}} \sum_{s=t}^T (T-s+1) = \frac{\tau(T-t+1)(T-t+2)}{2(T+1)^{3/2}},
\end{equation}
where we used $\sum_{s=t}^T T-s+1= \sum_{j=1}^{T-t+1} j$. The lower bound for $S_0$ follows by plugging in $t=0$ in the above bound and using a basic estimate.

For assertion (iii), the definition of $\tau_t$ and \eqref{eq: stbd1} give
\[\begin{aligned} 
    \sum_{t=0}^T \frac{\tau_t^2}{S_t} 
    = \frac{2\tau}{(T+1)^{3/2}} \sum_{t=0}^T \frac{T-t+1}{T-t+2} \leq  \frac{2\tau}{\sqrt{T+1}},
\end{aligned}\] 
where the last step is because the summands are upper bounded by $1$. This is the third estimate.

For assertion (iv), the definition of $\tau_t$ and \eqref{eq: stbd1} give
\[
\sum_{t=0}^{T-1} \frac{\tau_t^2}{S_{t+1}}
= \frac{2\tau}{(T+1)^{3/2}} \sum_{t=0}^{T-1} \frac{T-t+1}{T-t}
\leq \frac{4\tau T}{(T+1)^{3/2}}
\leq \frac{4\tau}{\sqrt{T+1}},
\]
since $\frac{T-t+1}{T-t} \leq 2$ for $t \leq T-1$.
\end{proof}

\begin{lemma}\label{proxboundsforstep4}
Let
$\tau_t=\frac{\tau}{\ln(t+2)\sqrt{t+2}}$
for some $\tau>0$. Then, for any $T\geq2$, we have
{\small
\[
\begin{alignedat}{2}
\text{(i)}\;\; \sum_{t=0}^T \tau_t^2
&\leq 3\tau^2,
\qquad &
\text{(ii)}\;\; S_0
&\geq \frac{\tau\sqrt{T+1}}{\sqrt{2}\ln(T+2)},
\\
\text{(iii)}\;\; \sum_{t=0}^T \frac{\tau_t^2}{S_t}
&\leq \frac{6\sqrt{2}\tau\ln(T+2)+33\tau}{\sqrt{T+1}},
\qquad &
\text{(iv)}\;\; \sum_{t=0}^{T-1} \frac{\tau_t^2}{S_{t+1}}
&\leq \frac{6\sqrt{2}\tau\ln(T+2)+32\tau}{\sqrt{T+1}}.
\end{alignedat}
\]
}
\end{lemma}
\begin{proof}
Assertion (i) follows from a standard estimation.
Now, we estimate $S_t$ (defined in the statement of \Cref{theorem: proxsgd_decreasingstepsize}) as
\[\begin{aligned}
S_t &= \sum_{s=t}^T \tau_s = \sum_{s=t}^T \frac{\tau}{\ln (s+2)\sqrt{s+2}} 
\geq \frac{\tau (T-t+1)}{\ln (T+2)\sqrt{T+2}}. 
\end{aligned}
\]

    Assertion (ii) follows by plugging in $t=0$ and using the fact that $\sqrt{T+2} \leq \sqrt{2} \sqrt{T+1}$.

    We next prove assertion (iv). In order to control $\sum_{t=0}^{T-1} \frac{\tau_t^2}{S_{t+1}}$, we will split the sum into two parts: $t=0,\ldots, \lfloor T/2 \rfloor - 1$ and $t=\lfloor T/2 \rfloor, \ldots, T-1$. Let $m = \lfloor T/2 \rfloor$. For the first half, the denominator $S_{t+1}$ can be bounded from below by $S_m$. In the second half, the numerator $\tau_t^2$ can be bounded by $\tau_m^2$. We will use these facts to bound the two parts separately.
    
    For the first part, we have
    \begin{equation}\label{eq: ubht5}
    \begin{aligned}
        \sum_{t=0}^{m - 1} \frac{\tau_t^2}{S_{t+1}} &\leq \frac{1}{S_m}\sum_{t=0}^{m - 1} \tau_t^2, \\
    \end{aligned}
    \end{equation}
    where we used the fact that $S_{t+1}$ is decreasing in $t$. Now, we estimate $S_m$ by lower bounding $\tau_s$ as 
    \[\begin{aligned}
        S_m &= \sum_{s=m}^T \tau_s \geq \frac{\tau(T - m + 1)}{\ln (T+2)\sqrt{T+2}} \geq \frac{\tau (T+1)}{2\ln (T+2)\sqrt{T+2}} \geq \frac{\tau \sqrt{T+1}}{2 \sqrt{2}\ln (T+2)},
    \end{aligned}\]
    where we used $T - m + 1 \geq (T+1)/2$ and $\sqrt{T+2} \leq \sqrt{2} \sqrt{T+1}$. Plugging this into \eqref{eq: ubht5} yields
    \begin{equation}\label{eq: first half}
        \sum_{t=0}^{m-1} \frac{\tau_t^2}{S_{t+1}} \leq \frac{2 \sqrt{2}\ln (T+2)}{\tau \sqrt{T+1}} \sum_{t=0}^{m-1} \tau_t^2 
        \leq \frac{6 \sqrt{2}\tau \ln (T+2)}{\sqrt{T+1}},
    \end{equation}
    where we use $m-1\leq T$ and assertion (i). 
    
    For the second half of the sum, we have $t\geq m$. First, note that since $\tau_t$ is non-increasing, we have that
    \begin{equation*}
    S_{t+1} = \sum_{i=t+1}^T \tau_i \geq (T-t) \tau_T.
    \end{equation*}
    With this, and $\tau_t$ being non-increasing, we have
    \begin{equation*}
    \sum_{t=m}^{T-1} \frac{\tau_t^2}{S_{t+1}} \leq \tau_m^2 \sum_{t=m}^{T-1} \frac{1}{S_{t+1}} \leq \tau_m^2 \sum_{t=m}^{T-1} \frac{1}{(T-t)\tau_T}.
    \end{equation*}
    Using now the definitions of $\tau_m, \tau_T$, and $\sum_{t=m}^{T-1} \frac{1}{T-t} \leq 1 + \ln(T+1)$, we obtain
    \begin{align}\label{eq: sum_tau22}
        \sum_{t=m}^{T-1} \frac{\tau_t^2}{S_{t+1}} &\leq  \frac{\tau \ln (T+2) \sqrt{T+2}}{(m+2)\ln^2 (m+2)} (1+\ln (T + 1)).
    \end{align}
First, it is trivial to obtain from $m = \lfloor T/2 \rfloor \geq T/2 - 1$ the estimate
    \begin{equation}\label{eq: estimate_m}
    m+2 \geq \frac{T+2}{2}.
    \end{equation}
    Second, it also follows, for $T \geq 2$, that
    \begin{equation}\label{eq: estimate_ln}
    \ln(m+2) \geq \frac{1}{2} \ln(T+2).
    \end{equation}
    Plugging \eqref{eq: estimate_m} and \eqref{eq: estimate_ln} back into \eqref{eq: sum_tau22} gives us
    \begin{equation}\label{eq: second half}
        \sum_{t=m}^{T-1} \frac{\tau_t^2}{S_{t+1}} \leq \frac{8 \tau (1+\ln (T + 1))}{\sqrt{T+2} \ln (T+2)} \leq \frac{32 \tau}{\sqrt{T+2}},
    \end{equation}
    where in the second inequality, we used $\frac{1 + \ln(T + 1)}{\ln(T+2)} \leq 4$ for $T \geq 0$. Combining \eqref{eq: first half} and \eqref{eq: second half} gives us
    \[
    \begin{aligned}
    \sum_{t=0}^{T-1} \frac{\tau_t^2}{S_{t+1}} &\leq \frac{6 \sqrt{2}\tau \ln (T+2) + 32 \tau}{\sqrt{T+1}}.\\
    \end{aligned}
    \]

    Finally, we prove assertion (iii). Since $S_t \geq S_{t+1}$ for $t=0,\ldots, T-1$ and $S_T = \tau_T$, splitting the sum and using assertion (iv), we have
    \[\begin{aligned}
        \sum_{t=0}^T \frac{\tau_t^2}{S_t} &= \sum_{t=0}^{T-1} \frac{\tau_t^2}{S_t} + \frac{\tau_T^2}{S_T} \leq \sum_{t=0}^{T-1} \frac{\tau_t^2}{S_{t+1}} + \tau_T \leq \frac{6 \sqrt{2}\tau \ln (T+2) + 32 \tau}{\sqrt{T+1}} + \frac{\tau}{\ln (T+2)\sqrt{T+2}}.
    \end{aligned}
    \]
    Using $1/\ln(T+2) < 1$ for $T \geq 1$ and $\frac{1}{\sqrt{T+2}} < \frac{1}{\sqrt{T+1}}$ for $T \geq 0$ gives the inequality.
\end{proof}
We have the following fact for recentering the \eqref{eq: BG} assumption. Its proof is trivial by using $\|\vx-\va\|^2 \leq 2\|\vx-\vb\|^2 + 2\|\vb-\va\|^2$.
\begin{fact}\label{fact: bg centering}
    If $\E\|\sgrad f(\vx)\|^2 \leq B^2 \|\vx-\va\|^2 + G^2$ for all $\vx$ and some fixed $\va$, then for any other point $\vb$, we have
    \[
        \E\|\sgrad f(\vx)\|^2 \leq 2B^2 \|\vx-\vb\|^2 + 2B^2 \|\vb-\va\|^2 + G^2.
    \]
\end{fact}

\vskip 0.2in
\bibliography{references}

\end{document}